\documentclass[10pt]{amsart}
\usepackage[margin=1.2in]{geometry}
\usepackage{colortbl}
\usepackage{color}
\usepackage{cite}

\newtheorem{theorem}{Theorem}[section]

\theoremstyle{definition}

\newtheorem{example}[theorem]{Example}

\theoremstyle{remark}
\newtheorem{remark}[theorem]{Remark}

\numberwithin{equation}{section}

\usepackage{graphicx}
\graphicspath{{./figures/}}
\usepackage{subcaption}
\usepackage{multirow}
\usepackage{bbm}
\usepackage{lipsum}
\usepackage{amsfonts}
\usepackage{graphicx}
\usepackage{epstopdf}
\usepackage{algorithmic}
\usepackage{subcaption}
\usepackage{amssymb}
\usepackage{mathtools}
\usepackage{bbm}

\newcommand{\minmod}{{\rm MM}}

\begin{document}
	
	\title{A local velocity grid conservative semi-Lagrangian schemes for BGK model}

	\author[S. Boscarino]{Sebastiano Boscarino}
	\address{Sebastiano Boscarino\\
		Department of Mathematics and Computer Science  \\
		University of Catania\\
		95125 Catania, Italy} \email{boscarino@dmi.unict.it}
	
	\author[S. Y. Cho]{Seung Yeon Cho}
	\address{Seung Yeon Cho\\
		Department of Mathematics and Computer Science  \\
		University of Catania\\
		95125 Catania, Italy}
	\email{chosy89@skku.edu}
	
	\author[G. Russo]{Giovanni Russo}
	\address{Giovanni Russo\\
		Department of Mathematics and Computer Science  \\
		University of Catania\\
		95125 Catania, Italy} \email{russo@dmi.unict.it}

	\maketitle
	
	\begin{abstract}
	
	Most numerical schemes proposed for solving BGK models for rarefied gas dynamics are based on the discrete velocity approximation. Since such approach uses fixed velocity grids, one must secure a sufficiently large domain with fine velocity grids to resolve the structure of distribution functions. When one treats high Mach number problems, the computational cost becomes prohibitively expensive. In this paper, we propose a velocity adaptation technique in the semi-Lagrangian framework for BGK model. The velocity grid will be set locally in time and space, according to mean velocity and temperature. We apply a weighted minimization approach to impose conservation. We presented several numerical tests that illustrate the effectiveness of our proposed scheme.
	\end{abstract}

	\section{Introduction}
	
%
	
	In the rarefied gas dynamics, the BGK model \cite{BGK} has been widely adopted as an approximation of the Boltzmann equation because of its simple structure. BGK collision operator is easier to compute than the Boltzmann one, and it allows efficient implementation of implicit schemes, therefore it can be used also when the Knudsen number is small. 
	Furthermore, using a penalty method, it allows the construction of efficient schemes for the full Boltzmann equation for small Knudsen number, as reported in \cite{BS,FJ}.
	
	Various numerical schemes have been proposed for solving the BGK model, based on the discrete velocity model (DVM). Among such schemes, our main interest concerns the class of semi-Lagrangian (SL) methods, which allow the use of large time step. In \cite{GRS, BCRY,CBRY1}, the combination of the semi-Lagrangian approach for the convection term is considered together with an implicit treatment of the collision term, and it enables us to avoid the CFL-type restriction while handling the stiffness problem coming from small Knudsen number. 

	Recently developed conservative semi-Lagrangian schemes allow accurate solutions on a wide range of Knudsen numbers, with very mild restrictions on the time step \cite{CBRY2}. A related convergence proof of the SL scheme can be found in \cite{RSY,RY,BCRY2}.
	
	In general, such approaches use fixed velocity grids, and one must secure a sufficient number of grid points in phase space to resolve the structure of the distribution function. When dealing with high Mach number problems, where large variation of mean velocity and temperature are present in the domain under consideration, the computational cost and memory allocation requirements become prohibitively large. To overcome such difficulty, local velocity grid methods have been developed  in the context of Eulerian based schemes \cite{BIP,BRULL201422}. We also refer to a recent work \cite{BP}, where the local velocity approach has been extended to a BGK model for gas mixtures.
	
	In this paper, we introduce a velocity adaption technique for the semi-Lagrangian scheme applied to the BGK model. The velocity grids will be set locally in time and space. We apply a weighted minimization approach to impose global conservation, generalizing the $L^2$-minimization technique introduced in \cite{gamba2009spectral}. We demonstrate the efficiency of the proposed scheme in several numerical examples.
	
	The outline of this paper is following. In Section \ref{sec: recon}, we review the conservative SL method for the BGK model \cite{CBRY1,CBRY2}. Then in Section \ref{sec weight}, we propose a weighted $L^2$-minimization approach to enforce conservation. Section \ref{sec 4} is devoted to the description of the local velocity grid approach in the semi-Lagrangian framework. Then, in Section \ref{sec cons poly}, we explain how we reconstruct numerical solutions for each cell. In Section \ref{sec second SL}, we describe a second order scheme. Finally, in Section \ref{sec: numeircal test}, we perform several numerical tests to demonstrate the efficiency of our methods.

\subsection{BGK model for the Boltmann equation}
For a small Knudsen number $\varepsilon \leq 10^{-4}$, the BGK model gives an good approximation for the Boltzmann equation because both equations lead to a same compressible Euler system in the limit $\varepsilon \rightarrow 0$. The BGK model replaces the collision term of the Boltzmann equation with a relaxation term. The BGK model is given by
\begin{equation}\label{bgk}
		\frac{\partial{f}}{\partial{t}} + v \cdot \nabla_x{f} = \frac{1}{\varepsilon}\left(\mathcal{M}(f)-f\right),
\end{equation}
where $f(x,v,t)$ denotes the number density of monatomic gas on a phase point  $(x,v) \in \mathbb{R}^{d_x} \times \mathbb{R}^{d_v}$ at time $0 \leq t \in \mathbb{R}_+$.
The local Maxwellian $\mathcal{M}(f)$ is given by 
	\begin{align*} 
		\mathcal{M}(f)(x,v,t):=\frac{\rho(x,t)}{\sqrt{\left(2 \pi RT(x,t) \right)^{d_v}}}\exp\left(-\frac{|v-U(x,t)|^2}{2RT(x,t)}\right).
	\end{align*}
where $R$ is the gas constant. Integration of the number density with respect to $\left(1,\,v,\,|v|^2/2,|v-U|^2/2\right)$ gives the information on macroscopic quantities such as mass $\rho(x,t)$, bulk velocity $U(x,t)$, total energy $E(x,t)$ and temperature $T(x,t)$ as follows:
\begin{align*}
		\begin{split}
			\bigg(\rho,\, \rho U,\, E,\,  \frac{\rho RT}{\gamma-1}\bigg)^\top& = \int_{\mathbb{R}^{d_v}} \left(1,\, v ,\,\frac{| v|^2}{2},\frac{|v-U|^2}{2}\right) f  dv,
		\end{split}
\end{align*}
Here we consider a monoatomic gas of unit mass for which the ratio of specific heats $\gamma$ is given by 
\[
\gamma= \frac{d_v+2}{d_v}.
\]
Note that the relaxation term still preserves fundamental properties of
the Boltzmann collision operator:
\begin{itemize}
	\item Collision invariance $ 1,v,|v|^2$:
		\begin{align*} 
			\int_{\mathbb{R}^{d_v}} (\mathcal{M}(f)-f)\begin{pmatrix}
				1\\v\\ |v|^2
			\end{pmatrix}dv=0.
	\end{align*}
	\item Conservation laws for mass, momentum and energy:
		\[
		\frac{d}{dt}  \int_{\mathbb{R}^{d_x} \times\mathbb{R}^{d_v} } f\begin{pmatrix}
			1\\v\\ |v|^2
		\end{pmatrix} \,dx\,dv=0,
		\]
	\item The H-theorem:
		\[
		\frac{d}{dt}  \int_{\mathbb{R}^{d_x}\times\mathbb{R}^{d_v}} f \log{\frac{1}{f}} \,dx\,dv   \geq 0.
		\]
	\item Taking integration of \eqref{bgk} with respect to $(1,v,|v|^2)$, we get 
	\begin{align*} 
		\begin{split}
			&\partial_t\langle f \rangle + \nabla \cdot \langle vf \rangle=0,\cr
			&\partial_t\langle vf \rangle + \nabla \cdot\langle v\otimes v f \rangle=0,\cr
			& \partial_t \bigg\langle \frac{|v|^2}{2}f \bigg\rangle + \nabla \cdot\bigg\langle\frac{|v|^2}{2}f\bigg\rangle =0, 
		\end{split}
	\end{align*} 	
	where $\langle g\rangle := \int_{\mathbb{R}^{d_v}} g(v) dv $. In the fluid regime $\varepsilon \rightarrow 0$, the solution $f$ tends to $\mathcal{M}$. Then, its macroscopic moments satisfy the compressible Euler system:
		\begin{align*}
			\begin{split}
				&\rho_t + \nabla \cdot(\rho u)=0,\cr
				&(\rho u)_t + \nabla \cdot(\rho u \otimes u+\rho RT I_d)=0,\cr
				& E_t + \nabla \cdot\left(\left(E+p\right)u\right) =0, 
			\end{split}
	\end{align*}
	where $p=\rho R T$ is pressure and $I_d$ is a $d_x\times d_x$ identity matrix. In the rest of this paper, we assume $d_x=d_v=1$.
\end{itemize}

\section{Review of the conservative SL scheme}\label{sec: recon}
\subsection{Conservative Reconstruction}\label{sec:cons_rec}
Here we briefly review the one-dimensional point-wise conservative reconstruction technique introduced in \cite{CBRY1,CBRY2}. Let us consider a uniform mesh size $\Delta x$ with grid points $x_i\equiv x_{min} + i \Delta x$, of the computational domain $[x_{min}, x_{max}]$. We denote by $\mathcal{I}$ the set of all space indices. Suppose that $u(x)= \frac{1}{\Delta x} \int_{x-\Delta x/2}^{x-\Delta x/2} \hat{u}(y)dy$ is the function we want to reconstruct. The procedure for the conservative reconstruction is given as follows:
\begin{enumerate}
	\item Given point-wise values $\{u_i\}_{i \in \mathcal{I}}$  for each $i \in \mathcal{I}$, we reconstruct a polynomial of even degree $k$:
	$$R_i(x)= \sum_{\ell=0}^k \frac{R_i^{(\ell)}}{\ell !}(x-x_i)^{\ell}$$ 
	which has the following properties:
	\begin{itemize}
		\item High order accurate in the approximation of smooth $\hat{u}(x)$ (see \cite{BCRY1}, Proposition 2.1):
		\begin{itemize}
			\item If $\ell$ is an even integer such that $0 \leq \ell\leq k$,
						\begin{align*}
					\begin{split}
						\hat{u}_i^{(\ell)}&=R_i^{(\ell)}+ \mathcal{O}(\Delta x^{k+2-\ell}). 
					\end{split}
			\end{align*}
			\item If $\ell$ is an odd integer such that $0\leq \ell < k$,
					\begin{align*}
					\begin{split}
						\hat{u}_i^{(\ell)}-\hat{u}_{i+1}^{(\ell)}&=R_i^{(\ell)}-R_{i+1}^{(\ell)}+ \mathcal{O}(\Delta x^{k+2-\ell}).
					\end{split}
			\end{align*}
		\end{itemize}
		\item Essentially non-oscillatory. 
		\item Positivity preserving.
		\item Conservative in the sense of cell averages:
				\[
			\frac{1}{\Delta x}\int_{x_{i-\frac{1}{2}}}^{x_{i+\frac{1}{2}}}R_i(x)\,dx = u_{i}.
			\] 
	\end{itemize}    
	\item Using the obtained values $R_i^{(\ell)}$ for $0 \leq \ell\leq k$, we approximate $u(x_{i+\theta})$, $\theta\in [0,1)$, to $O(\Delta x^{k+2})$ with
	\begin{align*}
			Q(x_i+\theta \Delta x)= \sum_{\ell=0}^k (\Delta x)^{\ell} \left(\alpha_{\ell}(\theta)R_{i}^{(\ell)} +\beta_\ell(\theta)R_{i+1}^{(\ell)}\right),
	\end{align*}
	where $\alpha_\ell(\theta)$ and $\beta_\ell(\theta)$ are given by 
		\begin{align*}
			\alpha_\ell(\theta) = \frac{1- (2\theta -1)^{\ell+1}}{2^{\ell+1}(\ell+1)!}, \quad 
			\beta_\ell(\theta) = \frac{(2\theta -1)^{\ell+1} - (-1)^{\ell+1} }  {2^{\ell+1}(\ell+1)!}.
	\end{align*}
	for $\theta \in [0,1)$.
\end{enumerate}
{In \cite{CBRY1}, we showed that CWENO polynomials  satisfy this conditions and we take it as basic reconstruction for the implementation of conservative SL schemes for BGK model in \cite{CBRY2}.

The technique can be easily extended to more dimensions, see \cite{CBRY1}.
}

\subsection{A semi-Lagrangian method for BGK model}
In a previous work \cite{CBRY2}, we introduced a  semi-Lagrangian method for BGK model, 
where we apply the conservative reconstruction in Section 
\ref{sec:cons_rec} to evaluate the distribution  function on off-grid points. In this section we review the first order SL method for the BGK model in \cite{CBRY2}. Applying implicit Euler method to its characteristic form, we get
\begin{align}\label{classical SL}
		f_{i,j}^{n+1}=\tilde{f}_{i,j}^n + \frac{\Delta t}{\varepsilon} \left(	\mathcal{M}_{ij}^{n+1}-f_{i,j}^{n+1}\right),
\end{align}
where $\tilde{f}_{i,j}^n$ is the approximation of $f(x_i - v_j \Delta t,v_j,t^n)$ which can be computed from  $\{f_{i,j}^n\}_{i \in \mathcal{I}}$ by a suitable reconstruction \cite{CBRY1,CBRY2} that enables us to preserve the global macroscopic moments. Here the local Maxwellian is computed by
\begin{align*}
			\mathcal{M}_{ij}^{n+1} =\frac{\rho_{i}^{n+1}}{\sqrt{\left(2 \pi T_{i}^{n+1} \right)^2}}\exp\left(-\frac{|v_{j}-U_{i}^{n+1}|^2}{2T_{i}^{n+1}}\right),
\end{align*}
with discrete macroscopic quantities:
\begin{align*}
		\left(\rho_{i}^{n+1},\, \rho_{i}^{n+1}U_{i}^{n+1},\, d_v \rho_{i}^{n+1} T_{i}^{n+1}\right) &:= \sum_{j \in \mathcal{J}} f_{i,j}^{n+1} \left(1,v_j,\big| v_j- U_{i}^{n+1} \big|^2\right)(\Delta v)^{d_v}.
\end{align*}
Thanks to the collision invariant $\phi_{j}:=(1,v_{j},\frac{|v_{j}|^2}{2})$,
one can compute \eqref{classical SL} explicitly. Multiplying the collision invariants to both sides of \eqref{classical SL} and taking summation over $j \in \mathcal{J}$, one obtains
\begin{align*} 
		\sum_{j \in \mathcal{J}} \big(f_{i,j}^{n+1} - \tilde{f}_{i,j}^n\big)\phi_{j} (\Delta v)^{d_v} = \frac{\Delta t}{\varepsilon}\sum_{j \in \mathcal{J}} \left(	\mathcal{M}_{ij}^{n+1}-f_{i,j}^{n+1}\right)\phi_{j} (\Delta v)^{d_v}.
\end{align*}
Here the right hand side can be negligible if the discrete summation is computed with a sufficiently refined grid on the appropriate velocity domain, because midpoint rule is spectrally accurate when applied to a Maxwellian.

Therefore, the discrete macroscopic quantities $\rho_i^n$, $U_i^n$ and $T_i^n$ can be replaced with
\begin{align*}
		\left(\rho_{i}^{n+1},\, \rho_{i}^{n+1}U_{i}^{n+1},\, E_{i}^{n+1}\right) &=\sum_{j \in \mathcal{J}} f_{i,j}^{n+1}\phi_{j} (\Delta v)^{d_v} = \sum_{j \in \mathcal{J}} \tilde{f}_{i,j}^n \phi_{j} (\Delta v)^{d_v}=:\left(\tilde{\rho}_{i}^{n+1},\, \tilde{\rho}_{i}^{n+1}\tilde{U}_{i}^{n+1},\, \tilde{E}_{i}^{n+1}\right),
\end{align*}
this further gives
\begin{align*}
		d_v \rho_{i}^{n+1}T_{i}^{n+1} &=\sum_{j \in \mathcal{J}} f_{i,j}^{n+1} \big|v_j-U_i^{n+1}\big|^2(\Delta v)^{d_v} = \sum_{j \in \mathcal{J}} \tilde{f}_{i,j}^n \big|v_j-\tilde{U}_i^{n+1}\big|^2 (\Delta v)^{d_v} = d_v\tilde{\rho}_{i}^{n+1}\tilde{T}_{i}^{n+1},
\end{align*}
Finally, one can update solution as follows:
\begin{align*} 
		f_{i,j}^{n+1}=\tilde{f}_{i,j}^n + \frac{\Delta t}{\varepsilon} \left(\tilde{\mathcal{M}}_{ij}^{n}-f_{i,j}^{n+1}\right),
\end{align*}
with
\begin{align*}
		\tilde{\mathcal{M}}_{ij}^{n} =\frac{\tilde{\rho}_{i}^{n}}{\sqrt{\left(2 \pi \tilde{T}_{i}^{n} \right)^{d_v}}}\exp\left(-\frac{|v_{j}-\tilde{U}_{i}^{n}|^2}{2\tilde{T}_{i}^{n}}\right).
\end{align*}

\section{Weighted $L^2$-minimization for moment correction}\label{sec weight}
When using semi-Lagrangian scheme for kinetic equations, the discrete conservation of the mass, momentum and energy may be lost. 
Lack of conservation has been analyzed in \cite{BCRY}, where it was found that the  non-linear weights used in the reconstructions break translation invariance causing lack of conservation in the distribution function, while computing the approximation of the moments by the discrete sums  destroys the exact conservation at the level of collision operator. In \cite{CBRY2}, we introduce a conservative reconstruction to solve the first problem, and adopt two techniques to maintain moment conservation: an approach based on Entropy minimization \cite{M} and one based on $L^2$ minimization \cite{gamba2009spectral}.

The first approach allows the construction of a conservative discrete Maxwellian, while $L^2$ minimization can be applied to more general distribution functions, so we adopt this approach in our paper. 
The technique is based on a constraint $L^2$-minimization, where an initial guess of the distribution function is slightly modified to impose conservation of the physical quantities. However, such a procedure may introduce negative values near the tails of the distribution.

Here we propose a
{\em weighted} $L^2$-minimization, which is more robust in preventing negative values of the distribution function. Let us consider reference mass, momentum and energy $\mathcal{U}:=(\rho,\rho U, E)^{\top} \in \mathbb{R}^{d_v+2}$.
Given an initial guess $f$, we consider a weight function $1/h$ and look for a solution $g$ of the following weighted $L^2$-minimization problem:
\begin{align}\label{min prb}
	\min_g \bigg\|f \circ \frac{1}{h} - g\bigg\|_2^2 \quad\text{ s.t }\quad Cg=\mathcal{U}
\end{align}
where $\circ$ denotes the componentwise multiplication and
$$f \equiv (f_1,f_2,\cdots,f_{(N_v+1)^{d_v}})^\top\in \mathbb{R}^{(N_v+1)^{d_v}}, \quad g =(g_1,g_2,\cdots,g_{(N_v+1)^{d_v}})^\top\in \mathbb{R}^{(N_v+1)^{d_v}}.$$
$$h \equiv (h_1,h_2,\cdots,h_{(N_v+1)^{d_v}})^\top\in \mathbb{R}^{(N_v+1)^{d_v}}.$$
$$\displaystyle C:=\begin{pmatrix}
	h_j(\Delta v)^{d_v}\\ h_j v_j (\Delta v)^{d_v} \\ h_j \frac{|v_j|^2}{2} (\Delta v)^{d_v}
\end{pmatrix}\in \mathbb{R}^{(d_v+2) \times (N_v+1)^{d_v}}, \quad \mathcal{U}= (\rho,\rho U, E)^{\top}\in \mathbb{R}^{(d_v+2) \times 1}.$$
Here $g$ is constructed as close as possible to the ratio of $f$ with respect to $h$ which corresponds to the macroscopic quantities $\mathcal{U}$, while $g \circ h$ gives the approximation of $f$ reproducing exactly the same discrete moments $\mathcal{U}$. Note that the use of weight $h_j\equiv 1$ for all $j$ leads to the classical $L^2$-minimization.

Now, we illustrate how we solve the weighted $L^2$-minimization problem \eqref{min prb}. As in \cite{gamba2009spectral} we use the method of Lagrange multiplier with the following Lagrangian $\mathcal{L}(g,\lambda)$:
$$\mathcal{L}(g,\lambda)= \left\|f\circ \frac{1}{h} - g\right\|_2^2 + \lambda^{\top}\left(C g-\mathcal{U}\right).$$
We first find the stationary points of  $\mathcal{L}$:
\begin{align*}
	\begin{array}{rlrl}
		\nabla_g \mathcal{L} =0 &\Leftrightarrow& g&= f\circ \frac{1}{h}+ \frac{1}{2}C^{\top} \lambda\cr
		\nabla_\lambda \mathcal{L} = 0 &\Leftrightarrow& C g&=\mathcal{U}.
	\end{array}
\end{align*}
From these, we explicitly compute $\lambda$ as follows:
\begin{align*}
	\lambda= 2(CC^{\top})^{-1}\left(\mathcal{U}-C\left(f\circ \frac{1}{h}\right)\right).
\end{align*}
Here the matrix $CC^{\top}$ is invertible because it is symmetric and positive definite. Consequently,
\begin{align*} 
	g\circ h= f+C^{\top}(CC^{\top})^{-1}\left(\mathcal{U}-C\left(f\circ \frac{1}{h}\right)\right)\circ h.
\end{align*}

The additional computational cost for this procedure with respect to the classical approach \cite{gamba2009spectral} is just a few component-wise division and multiplication of vectors of size $(N_v+1)^{d_v}$.

%

\begin{figure}[t]
	\centering
	\begin{subfigure}[t]{0.45\linewidth}
		\includegraphics[width=1\linewidth]{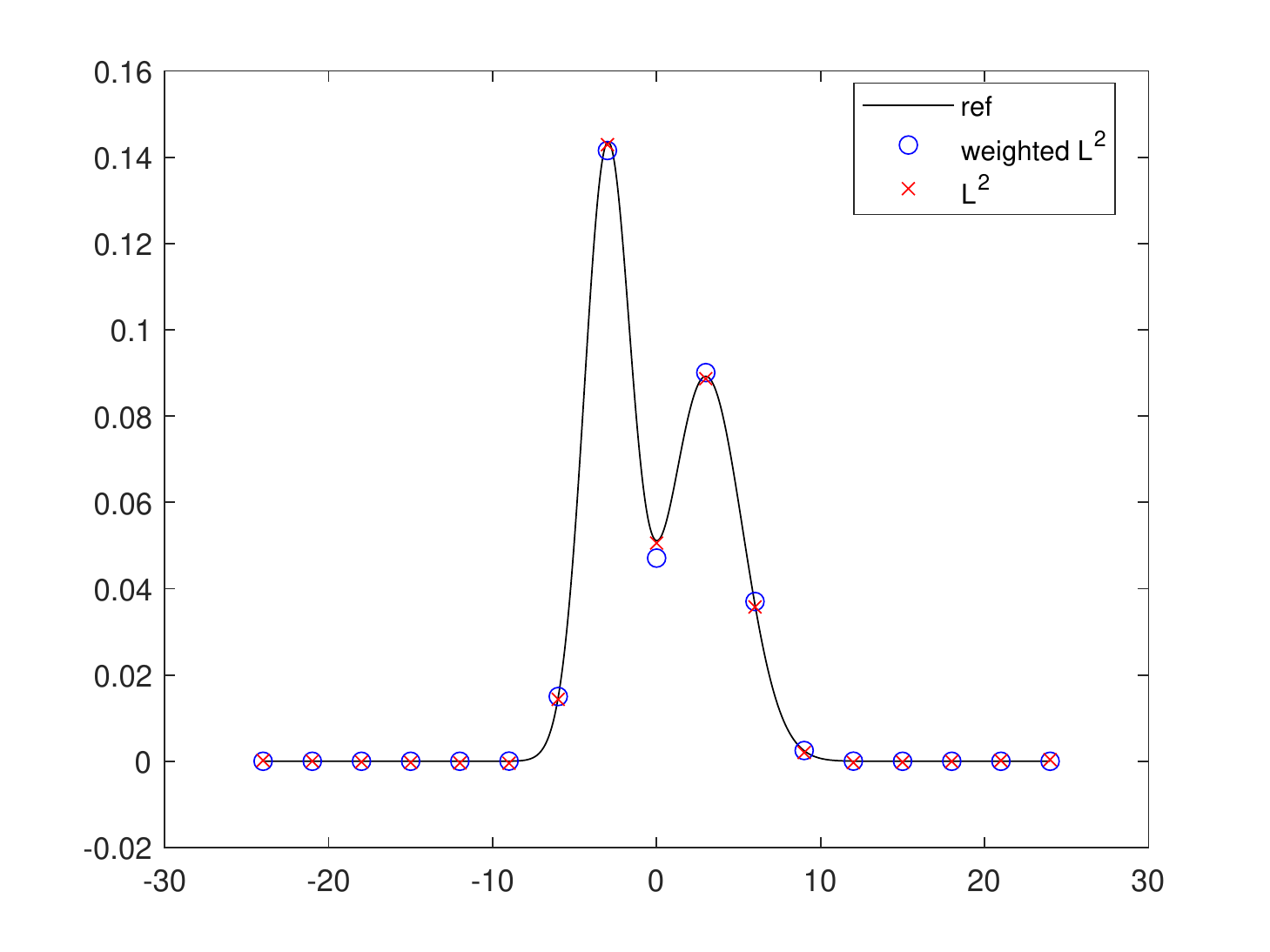}
	\end{subfigure}	
	\begin{subfigure}[t]{0.45\linewidth}
		\includegraphics[width=1\linewidth]{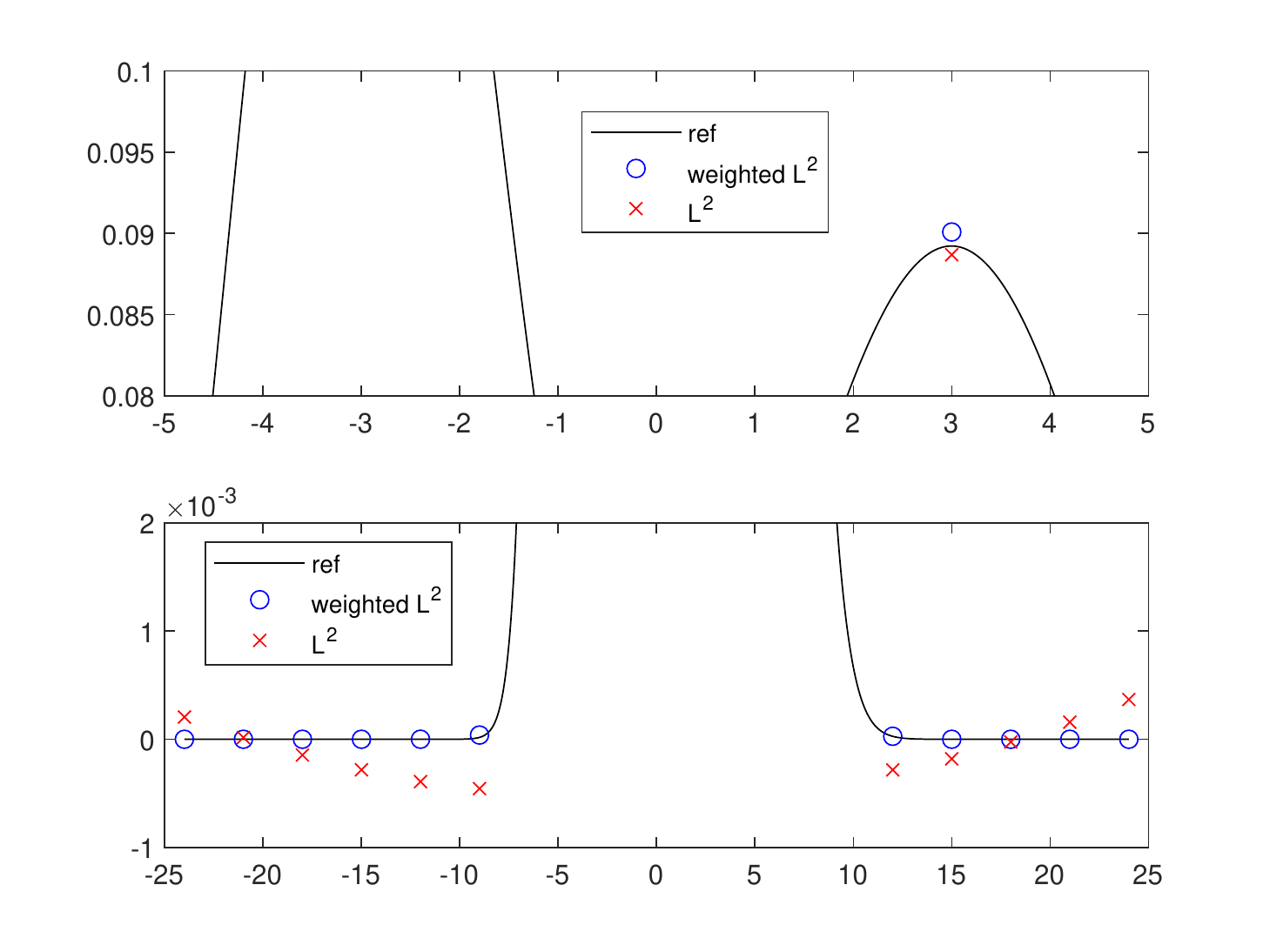}
	\end{subfigure}	
	\caption{Comparison of weighted and non-weighted $L^2$-minimization.}\label{comparison of weight}
\end{figure}

As an application of this approach, we present an example in Figure \ref{comparison of weight},
where we compare the results of weighted and non-weighted $L^2$-minimization techniques.
\begin{example}
	As a reference solution, we consider a distribution function:
	\begin{align*}
		f(v)&= \mathcal{M}_1(v)+\mathcal{M}_2(v)
	\end{align*}
	where
	\begin{align*}
		\mathcal{M}_1(v)&=\frac{0.5}{\sqrt{4 \pi}}\exp\left(-\frac{|v+3|^2}{4}\right),\quad \mathcal{M}_2(v)= \frac{0.5}{\sqrt{10 \pi}}\exp\left(-\frac{|v-3|^2}{10}\right).
	\end{align*}
	Then, the macroscopic moments of the distribution function $f$ are given by {$\rho=1$, $\rho U=0, E=6.25$}. Assume that we are given $f(v_j)$ on grid points $v_j \in \mathcal{V}=\{-24,-21,...,21,24\}$. For the weighted $L^2$-minimization approach, we use the Maxwellian constructed by $\rho, U, E$.
	In Figure \ref{comparison of weight}, we compare two solutions with the reference solution. Both solutions preserve the reference macroscopic quantities. As expected, the weighted approach tends to prevent the negative values on the tail of the distribution while it gives slightly less accurate solutions near the peaks of the distribution function.
\end{example}


\section{Semi-Lagrangian method based on local velocity grid approach}\label{sec 4}
In this section, we present semi-Lagrangian method based on the use of local velocity grid approach. Before proceeding we introduce the notation we will use throughout this paper. Let us consider a fixed time step $\Delta t$ and denote the $n$th time step by $t^n:=n \Delta t$. For space, we assume a uniform mesh size $\Delta x$ and each node is defined by $x_i=x_{min}+\left(i-\frac{1}{2}\right)\Delta x$. For each $i$ and $n$, let us denote the set of velocity grid points by $\mathcal{G}_i^{n}$ with uniform mesh size $\Delta v_i^n$. The set of indices corresponding to space and to each $\mathcal{G}_i^{n}$ will be denoted by $\mathcal{I}$ and $\mathcal{J}_i^{n}$, respectively. A cell assigned to each $(i,j) \in \mathcal{I} \times \mathcal{J}_i^n$ will be defined by  $I_{ij}^n:=\left[x_{i}-\Delta x/2, x_{i}+\Delta x/2\right] \times \left[v_{j}-\Delta v_i^n/2, v_{j}+\Delta v_i^n/2\right]$. The notation $\tilde{I}_{ij}^n$ is the region such that $$\tilde{I}_{ij}^n:=\{(x,v) \in \mathbb{R}^2 |\,  x_{i-\frac{1}{2}}- v \Delta t\leq x \leq  x_{i+\frac{1}{2}} - v \Delta t,\quad   v_{i,j}^n-\Delta v_i^n/2\leq v \leq v_{i,j}^n+\Delta v_i^n/2.
\}.$$

Now, we describe how to apply local velocity
adaptation approach to a first order SL scheme.

\subsection{First order semi-Lagrangian method for BGK model}\label{sec first SL}
Here we consider one-dimension in space $d_x=1$ and velocity $d_v=1$. To apply semi-Lagrangian scheme to \eqref{bgk}, we consider its characteristic form:
\begin{equation}\label{SL}
	\begin{array}{l}
		\displaystyle \frac{df}{ds} =\frac{1}{\varepsilon} \left(\mathcal{M}(f) - f\right), \quad t >0, \quad x,v \in \mathbb{R} \\[3mm] 
		\displaystyle \frac{dX}{ds} = v, \quad X(t) = x, \quad X(s) = x-v(t-s).
		\end{array}
\end{equation}

For treatment of the stiffness coming from $\varepsilon$, we begin by applying the implicit Euler scheme to \eqref{SL}. Then, we advance the solution from $t^n$ to $t^{n+1}$ using the following discretization:
\begin{align*}
	f^{n+1}(x,v)=f^n(x-v\Delta t,v) + \frac{\Delta t}{\varepsilon} \left(\mathcal{M}^{n+1}(x,v)-f^{n+1}(x,v)\right),
\end{align*}
where $f^n$ denotes the solution at time  $t^n$.


Our scheme can be explained by five steps. \textbf{Step 1}, we first predict $\rho_i^*$, $U_i^*$, $E_i^*$ and $T_i^*$ for each cell using the grid points $\mathcal{G}_i^n$ previously defined. 
\textbf{Step 2} we use the mean velocity $U_i^*$ and temperature $T_i^*$ to define a new grid points $\mathcal{G}_i^{n+1}$ for each $i$. 
\textbf{Step 3}, we improve the computation of $\rho_i^{n+1}$, $U_i^{n+1}$, $E_i^{n+1}$ and $T_i^{n+1}$ for each cell $i$. 
\textbf{Step 4}, we correct the Maxwellians
and the numerical solution using the weighted $L^2$-minimization technique to impose the correct conservation.
\textbf{Step 5} we update the numerical solution.

\begin{figure}[t]
	\centering
	\begin{subfigure}[t]{0.8\linewidth}
		\includegraphics[width=1\linewidth]{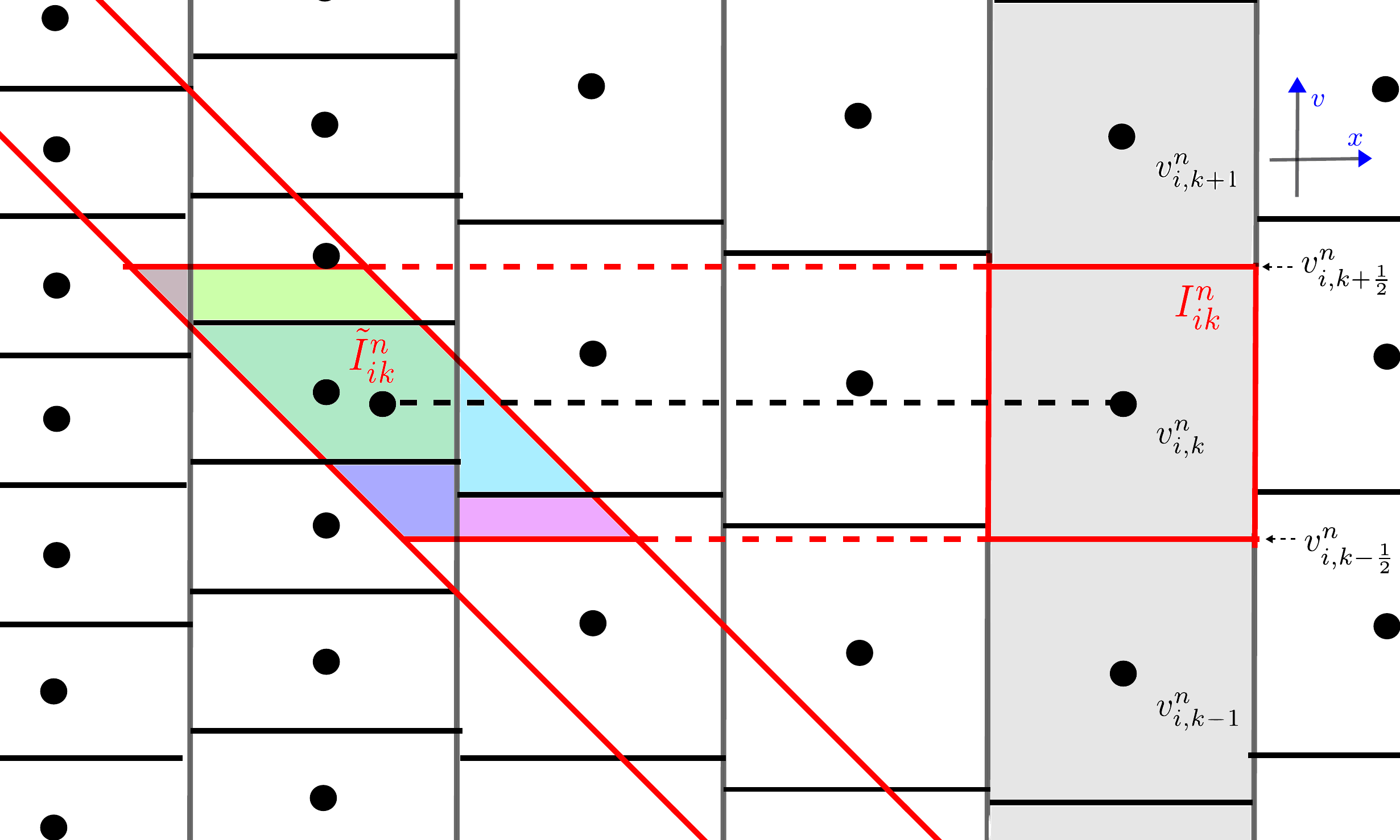}
	\end{subfigure}	
	\caption{Phase space discretization at time $t=t^{n}$ and computation of moment prediction by integration over the cell $\tilde{I}_{i,k}^n$ (\textbf{Step 1}). In each coloured patch, we use three different polynomials $P^n$, $Q^n$, $R^n$ for the computation of the contribution to $\rho_i^*$, $\rho_i^*U_i^*$, $E_i^*$. The grey region is the d }\label{parallelogram}
\end{figure}

\noindent{\textbf{Step 1: Prediction of macroscopic quantities at $t^{n+1}$ using $\mathcal{G}_i^n$}.} 
We first need to define the new velocity nodes for time $t^{n+1}$. For this, we begin by applying the implicit Euler method to \eqref{SL}:
\begin{align*}
	f^{n+1}(x,v)=f^n(x-v\Delta t,v) + \frac{\Delta t}{\varepsilon} \left(\mathcal{M}^{n+1}(x,v)-f^{n+1}(x,v)\right),
\end{align*}
where $f^n$ denotes the solution at time $t^n$. Taking integration over $I_{i}\times \mathbb{R}$, we get
\begin{align*} 
	\int_{\mathbb{R}}\int_{I_{i}}f^{n+1}(x,v)\begin{pmatrix}
		1  \\ v \\ v^2
	\end{pmatrix}\,dx\,dv&=\int_{\mathbb{R}}\int_{I_{i}}f^n(x-v\Delta t,v)\begin{pmatrix}
	1  \\ v \\ v^2
\end{pmatrix}\,dx\,dv=\int_{\mathbb{R}}\int_{\tilde{I}_{i}(v,\Delta t)}f^n(x,v)\begin{pmatrix}
	1  \\ v \\ v^2
\end{pmatrix}\,dx\,dv
\end{align*}
where $$I_i=\left[x_{i-\frac{1}{2}}, x_{i+\frac{1}{2}}\right], \quad \tilde{I}_{i}(v,\Delta t)=\left[x_{i-\frac{1}{2}}-v\Delta t, x_{i+\frac{1}{2}}-v\Delta t\right].$$
Based on this, we precompute $\rho_i^{*}$, $u_i^{*}$ and $E_i^{*}$ and temperature $T_i^{*}$:
\begin{align}\label{com rhouE}
\begin{split}
	\begin{pmatrix}
		\rho_i^{*}\\ \rho_i^{*}U_i^{*}\\ 2E_i^{*}
	\end{pmatrix}&:=\frac{1}{\Delta x}\int_{\mathbb{R}}\int_{\tilde{I}_{i}(v,\Delta t)}f^n(x,v)\begin{pmatrix}
		1  \\ v \\ v^2
	\end{pmatrix}\,dx\,dv \approx \sum_{k\in \mathcal{J}_i^n}\frac{1}{\Delta x}\int_{\tilde{I}_{i,k}^n}\begin{pmatrix}
		P^n(x,v)  \\ Q^n(x,v) \\ R^n(x,v)
	\end{pmatrix}\,dx\,dv,\cr
	d_vRT_i^{*}&:=(2E_i^{*}- \rho_i^{*}U_i^{*})/\rho_i^{*}.
\end{split}
\end{align}
This relation means that we compute macroscopic quantities based on the projected area of the grey region along characteristics. The parallelogram shaped domain $\tilde{I}_{i,k}^n$ is the projected area of $I_{ik}^n$ along the characteristics (see Figure \ref{parallelogram}). To construct a first order scheme, it is enough to take
\begin{align}\label{PQR piecewie}
	P^n(x,v)\equiv \sum_{i,k\in \mathcal{J}_i^n}f_{ik}^n \chi_{ik}^n(x,v),\quad Q^n(x,v)\equiv \sum_{i,k\in \mathcal{J}_i^n}v_kf_{ik}^n \chi_{ik}^n(x,v),\quad R^n(x,v)\equiv \sum_{i,k\in \mathcal{J}_i^n}v_k^2f_{ik}^n \chi_{ik}^n(x,v).
\end{align}
where $\chi_{ik}^n(x,v)$ is the characteristic function on $I_{ik}^n$. The integration in \eqref{com rhouE} will be obtained by summation of the integrals for trapzoidal patches. An explicit formula of such integral is illustrated in the Appendix \eqref{trap exact integral}.

For high order schemes, in section \ref{sec parallel}-\ref{sec con poly} we explain how we construct polynomials $P^n(x,v)$, $Q^n(x,v)$ and $R^n(x,v)$. 
%

\noindent{\textbf{Step 2: Choice of local velocity grids $\mathcal{G}_i^{n+1}$.}
Since our interests lie on the problems for small Knudsen number $\varepsilon$, we consider the situation when the shape of the distribution function is close to a local Maxwellian. Since it takes the form of the normal distribution whose mean is $U$ and standard deviation is $\sqrt{RT}$, most of the distribution function is concentrated in the interval
\[
U_i^* - \alpha\sqrt{RT_i^*} \leq v  \leq U_i^* + \alpha\sqrt{RT_i^*}.
\]
Here, a sufficiently large number $\alpha$ guarantees the approximate conservation of mass/momentum/energy. The $\alpha$ will be chosen large enough so that distributions outside this interval is acceptably small and can be negligibly small. To define velocity grids $\mathcal{G}_i^{n+1}$, we also need to set a size of mesh $\Delta v_i^{n+1}$. 

To resolve the shape of the Maxwellian, we impose that
between the two inflexion points of the Maxwellian there should be at least three grid points:
\begin{align}\label{dv cond}
	\Delta v_i^{n+1} \leq \sqrt{RT_i^*}.
\end{align}
Now, based on this, we compute
$N_v=2 \left\lceil \alpha\sqrt{RT_i^*}\big/\Delta v_i^{n+1}\right\rceil$
\begin{align*}
	v_{min,i}^{n+1}:= U_i^* -\frac{N_v}{2} \Delta v_i^{n+1},\quad v_{max,i}^{n+1}=U_i^* +\frac{N_v}{2} \Delta v_i^{n+1},
\end{align*}	
Here, for any $\xi \in \mathbb{R}$, $\lceil{\xi}\rceil$ means the smallest integer equal or greater than $\xi$. 
Finally, we newly define a set of grid points $\mathcal{G}_i^{n+1}$ by
\begin{align*}
	\mathcal{G}_i^{n+1}:= \{v_{i,j}^{n+1}|v_{i,j}^{n+1}:=v_{min,i}^{n+1} + j\Delta v_i^{n+1}, 0 \leq j\leq N_v\}.
\end{align*}

This, however, may be insufficient to resolve solutions where abrupt change in velocity and temperature appear.
Moreover, since we are considering a semi-Lagrangian framework, the information of nearest cells should be involved considering the CFL number. For this reason, we use
\[
T_{i,max}^* = \max_{-\delta\leq j\leq \delta}{T_{i+j}^*}, \quad T_{i,min}^* = \min_{-\delta\leq j\leq \delta}{T_{i+j}^*}, \quad 	\Delta v_i^{n+1} \leq \sqrt{RT_{i,min}^*}.
\]
where $\delta=\lceil{\text{CFL}}\rceil+1$. The additional number $1$ is due to the use of the integral based on parallelogram.

We determine $\Delta v_i^{n+1}$ as follows:
\begin{align*} 
	\Delta v_i^{n+1} = \beta \sqrt{RT_{i,min}^*},
\end{align*}
which defines the local resolution of the grid in velocity. For numerical simulation, we set $\beta=0.5$.
\begin{remark}
		When the shape of the distribution function is far from a local Maxwellian, its numerical support may be  relatively large, that is, the bounds could be far from the bounds computed from the Maxwellian.	For treating this problem, in \cite{BRULL201422,BP} the authors did as follows: once a transport step is performed for the prediction of the distribution function on the boundary of newly defined velocity domain, its relative scale is compared with a given tolerance. If the value is larger than the tolerance, the new velocity node is considered and the procedure is repeated until the values on endpoints become small enough. In our approach, we could treat the problem similarly. In the following step, we predict $f_{i,j}^{n+1}$ for each newly defined velocity node $v_{i,j}^{n+1}$, $j=0,...,N_v$. Then, additionally, we compute $f_{i,-1}^{n+1}$ (or $f_{i,N_v+1}^{n+1}$). If the value of $f_{i,-1}^{n+1}$ (or $f_{i,N_v+1}^{n+1}$) is larger than a tolerance $tol$, i.e.
		$$\frac{f_{i,-1}^{n+1}}{\max_{j\in \mathcal{J}_i^{n+1} }{f_{i,j}^{n+1}}} > tol \quad (\text{or} \, \frac{f_{i,N_v+1}^{n+1}}{\max_{j\in \mathcal{J}_i^{n+1} }{f_{i,j}^{n+1}}}>tol),$$
		we include the velocity node $j=0$ (or $j=N_v+1$) in $G_{i}^{n+1}$ and repeat the procedure until the value of distribution becomes small near boundary. As an alternative approach, we also refer to the work \cite{BIP}. 
\end{remark}

\begin{remark}
In this paper we shall use a constant time step $\Delta t$. This means that the CFL number, defined as 
\[
	{\rm CFL} = \frac{V_{\rm max}\Delta t}{\Delta x}
\]
will in general be not uniform, since the maximum grid velocity is different for different cells.
\end{remark}

\noindent{\textbf{Step 3: Correction of macroscopic quantities at $t^{n+1}$ and computation of $\tilde{f}_{ij}^n$ and $\mathcal{M}_{ij}^{n+1}$.}
We improve the prediction of the moments at time $t^{n+1}$ as follows. We first compute $\tilde{f}_{i,j}^n$ as an approximation of $f(x_i - v_j \Delta t,v_j,t^n)$:
\begin{align}\label{tilde f}
	\begin{split}
		\tilde{f}_{i,j}^n=\frac{1}{\Delta x \Delta v_i^{n+1}}\int_{\tilde{I}_{i,j}^{n+1}}P^n(x,v) \,dx\,dv.
	\end{split}
\end{align}
That is, we compute $\tilde{f}_{i,j}^n$ as the average of $P^n$ on $\tilde{I}_{i,j}^{n+1}$ (see Figure \ref{parallelogram np1}).

Now, we use this to compute $\rho_i^{n+1}$
\begin{align}\label{com rho}
	\begin{split}
		\rho_i^{n+1}&:=\sum_{j\in \mathcal{J}_i^{n+1}}\tilde{f}_{i,j}^n \Delta v_i^{n+1}.
	\end{split}
\end{align}
The other quantities $U_i^{n+1}$, $E_i^{n+1}$ and temperature $T_i^{n+1}$ are obtained by
\begin{align}\label{com rhouE np1}
	\begin{split}
		\begin{pmatrix}
			 \rho_i^{n+1}U_i^{n+1}\\ 2E_i^{n+1}
		\end{pmatrix}&:=\sum_{j\in \mathcal{J}_i^{n+1}}\frac{1}{\Delta x}\int_{\tilde{I}_{i,j}^{n+1}}\begin{pmatrix}
			 Q^n(x,v) \\ R^n(x,v)
		\end{pmatrix}\,dx\,dv,\cr
		d_vRT_i^{n+1}&:=(2E_i^{n+1}- \rho_i^{n+1}U_i^{n+1})/\rho_i^{n+1}.
	\end{split}
\end{align}
This step is needed because the new velocity domain (grey region in Figure \ref{parallelogram np1}) could be quite different. Notice that this step is not expensive
because the polynomials $P^n$, $Q^n$ and $R^n$ are already known from \textbf{Step 1}.
Then, the local Maxwellian $\mathcal{M}_{ij}^{n+1}$ is computed by
\begin{align*} 
	\mathcal{M}_{ij}^{n+1} =\frac{\rho_{i}^{n+1}}{\sqrt{\left(2 \pi T_{i}^{n+1} \right)^2}}\exp\left(-\frac{|v_{j}-U_{i}^{n+1}|^2}{2T_{i}^{n+1}}\right).
\end{align*}

\begin{figure}[t]
	\centering
	\begin{subfigure}[t]{0.8\linewidth}
		\includegraphics[width=1\linewidth]{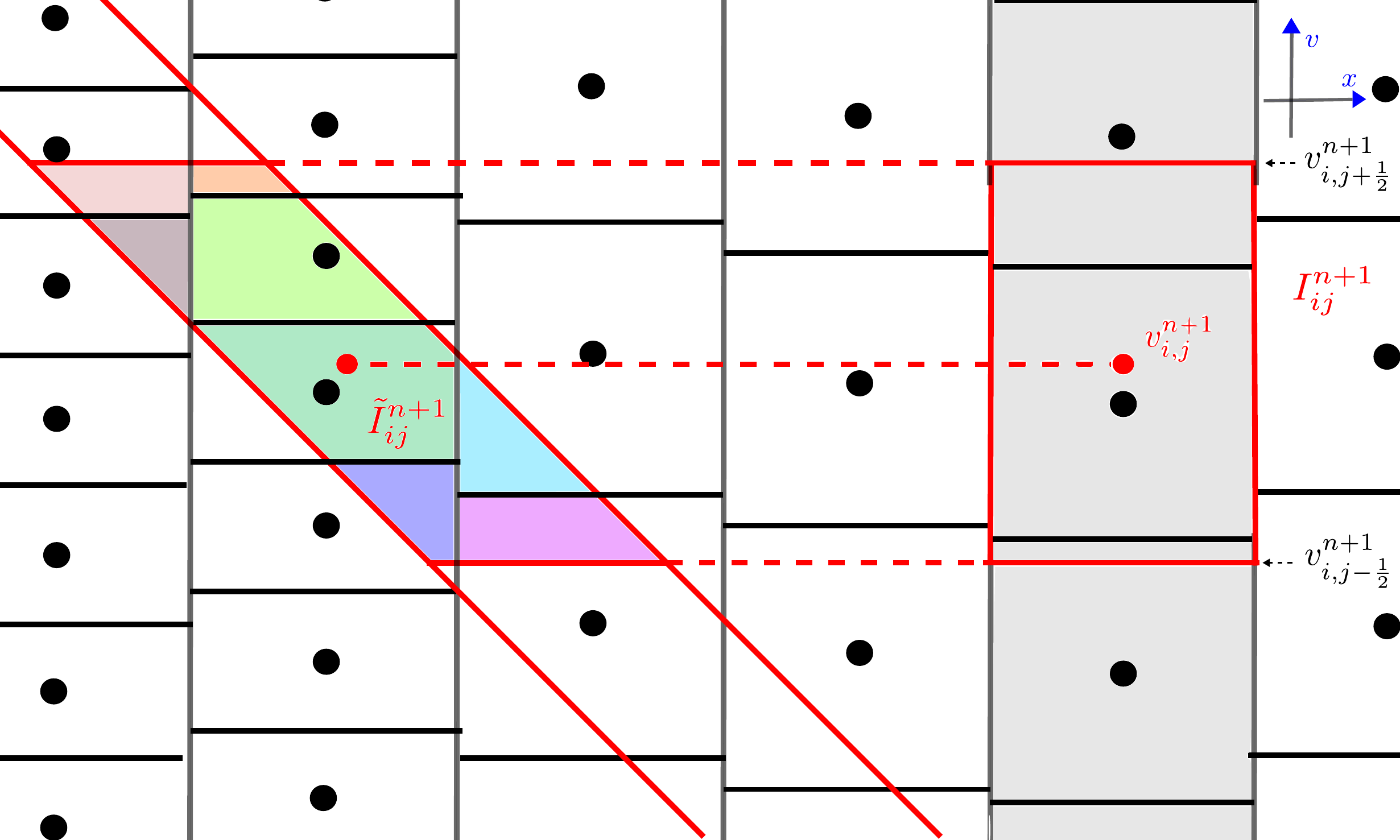}
	\end{subfigure}	
	\caption{Computation of the moments on the new velocity grid obtained by integration over the region $\tilde{I}_{i,j}^{n+1}$ (\textbf{Step 3}). In each patch, we use three different polynomials $P^n$, $Q^n$, $R^n$ for the computation of the contribution to $\rho_i^{n+1}$, $\rho_i^{n+1}U_i^{n+1}$, $E_i^{n+1}$.
}\label{parallelogram np1}
\end{figure}

\noindent{\textbf{Step 4: Correction of $\tilde{f}_{ij}^{n}$ and $\mathcal{M}_{ij}^{n+1}$ using a weighted $L^2$-minimization.}
Now, for each newly defined $v_j\in \mathcal{G}_i^{n+1}$ we attempt to compute solution as
\begin{align}\label{2D IE}
	f_{i,j}^{n+1}=\tilde{f}_{i,j}^n + \frac{\Delta t}{\varepsilon} \left(\mathcal{M}_{ij}^{n+1}-f_{i,j}^{n+1}\right).
\end{align}

The relation \eqref{2D IE} implies that the discrete summation of $f_{ij}^{n+1}$ and $\tilde{f}_{ij}^{n}$ over $1,v_j,v_j^2$ should reproduce the same macroscopic quantities. However, this may not be a good approximation if there are not enough grid points in velocity space.

In order to enforce conservation, we use the weighted $L^2$-minimization technique described in Section \ref{sec weight}.
First, we solve 
\begin{align}\label{min prb sl1}
	\min_p \bigg\|  {\bf{1}} - p_i^{n+1}\bigg\|_2^2 \quad\text{ s.t }\quad C_i^{n+1}p_i^{n+1}=\begin{pmatrix}
		\rho_i^{n+1}\\\rho_i^{n+1}U_i^{n+1}\\E_i^{n+1}
	\end{pmatrix}
\end{align}
where
$${\bf{1}}\equiv(1,1,\cdots,1)^\top\in \mathbb{R}^{(N_v+1)^{d_v}}.$$
$$p_i^{n+1} \equiv(p_{i,1}^{n+1},p_{i,2}^{n+1},\cdots,p_{i,N_v+1}^{n+1})^\top\in \mathbb{R}^{(N_v+1)^{d_v}}.$$
$$\displaystyle C_i^{n+1}:=\begin{pmatrix}
	\mathcal{M}_{i,j}^{n+1}(\Delta v)^{d_v}\\ \mathcal{M}_{i,j}^{n+1} v_j (\Delta v)^{d_v} \\ \mathcal{M}_{i,j}^{n+1} \frac{|v_j|^2}{2} (\Delta v)^{d_v}
\end{pmatrix}\in \mathbb{R}^{(d_v+2) \times (N_v+1)^{d_v}}.$$
where $v_j \in \mathcal{G}_i^{n+1}$ and $|\mathcal{G}_i^{n+1}|=N_v+1$. 
Then, we replace $\{\mathcal{M}_{i,j}^{n+1}\}_{j\in \mathcal{J}_i^{n+1}}$ with $M_i^{n+1}:=p_i^{n+1} \circ \mathcal{M}_{i}^{n+1}$.

Furthermore, notice that $\tilde{f}_{i,j}^{n}$ can be used only for computing the density $\rho_i^{n+1}$ by \eqref{tilde f}. In order to preserve conservation, we use again the weighted $L^2$-minimization and replace $\tilde{f}_{i,j}^{n}$ by a new distribution function $g_{i,j}^n$ which is the closest one to $\tilde{f}_{i,j}^{n}$ that provides correct moments $\rho_i^{n+1}$, $\rho_i^{n+1}U_i^{n+1}$ and $E_i^{n+1}$
previously computed in \eqref{com rho} and \eqref{com rhouE np1}.
Then, we find $g_{i}^n := q_{i}^{n+1} \circ \tilde{f}_{i}^{n}$
by solving
\begin{align*} 
	\min_q \bigg\|\tilde{f}_{i}^{n} \circ \frac{1}{M_i^{n+1}} - q_i^{n+1}\bigg\|_2^2 \quad\text{ s.t }\quad D_i^{n+1}q_i^{n+1}=\begin{pmatrix}
		\rho_i^{n+1}\\\rho_i^{n+1}U_i^{n+1}\\E_i^{n+1}
	\end{pmatrix}
\end{align*}
where
$$q_i^{n+1} \equiv(q_{i,1}^{n+1},q_{i,2}^{n+1},\cdots,q_{i,N_v+1}^{n+1})^\top\in \mathbb{R}^{(N_v+1)^{d_v}}.$$
$$\displaystyle D_i^{n+1}:=\begin{pmatrix}
	M_i^{n+1}(\Delta v)^{d_v}\\ M_i^{n+1} v_j (\Delta v)^{d_v} \\ M_i^{n+1} \frac{|v_j|^2}{2} (\Delta v)^{d_v}
\end{pmatrix}\in \mathbb{R}^{(d_v+2) \times (N_v+1)^{d_v}}.$$
where $M_i^{n+1}$ is the solution obtained from \eqref{min prb sl1}.

\noindent{\textbf{Step 5: Update of solution $f_{ij}^{n+1}$}.}
The final step is to update solution using 
\begin{align*} 
	f_{i,j}^{n+1}=g_{i,j}^{n}  + \frac{\Delta t}{\varepsilon} \left(M_{ij}^{n+1}-f_{i,j}^{n+1}\right).
\end{align*}

\section{Conservative treatment of the transport.}\label{sec cons poly}
In this section, we describe how to treat the transport part of the equation \eqref{bgk} in the framework of semi-Lagrangian scheme with velocity adaptation. The technique allows us to preserve the discrete moments such as mass/momentum/energy. For simplicity, we assume periodic boundary conditions. 
For standard conservative SL schemes based on global grids, such as the ones illustrated in Section \ref{sec: recon}, we use the same velocity grids $\mathcal{G}_i^n=\mathcal{G}$ with same mesh size $\Delta v_i^n=\Delta v$ for each spatial cell $i$ and time step $t^n$.
In this case, the conservative reconstruction technique of Section \ref{sec:cons_rec}
satisfies
\begin{align}\label{local summation}
	\sum_{i}f^{n+1}_{ij}\phi_j\Delta v_i^{n+1}=\sum_{i}\tilde{f}^{n}_{ij}\phi_j\Delta v_i^n 
	=\sum_{i}f^{n}_{ij}\phi_j\Delta v_i^n ,\quad \forall v_j\in \mathcal{G}
\end{align}
where $\phi_j=1,v_j,v_j^2/2$ and hence global conservation holds:
\begin{align}\label{global summation}
	\sum_i\begin{pmatrix}
		\rho^{n+1}_{i}\\
		m^{n+1}_{i}\\
		E^{n+1}_{i}
	\end{pmatrix}\Delta x=\sum_{i,j}\tilde{f}^{n}_{ij}\phi_j\Delta v \Delta x
	=\sum_{i,j}f^{n}_{ij}\phi_j\Delta v \Delta x.
\end{align}
The reason conservation is guaranteed is due to the fact that in all space cells there are the same discrete velocities, and interpolation is needed only on space. Once the conservative reconstruction is adopted, conservation is automatically guaranteed. 

However, relation \eqref{local summation} may not hold if one uses local velocity grid approach, because each space cell will have a different set of velocities, and therefore interolation at the level of velocity is needed, so if this is not properly done global conservation may be lost. In the rest of this section, we show how to restore global conservation
with a suitable choice of piecewise polynomials.
\subsection{Approximation based on a parallelogram}\label{sec parallel}
In order to ensure the global conservation \eqref{global summation}, an alternative way is to compute macroscopic quantities with \eqref{com rho} and \eqref{com rhouE np1}. Here we consider a more general form of \eqref{PQR piecewie}:
$$\begin{pmatrix}
	P^n(x,v)\\Q^n(x,v)\\R^n(x,v)
\end{pmatrix}=\sum_{ik} \begin{pmatrix}
P_{ik}^{n}(x,v)\\Q_{ik}^{n}(x,v)\\R_{ik}^{n}(x,v)
\end{pmatrix} \chi_{ik}(x,v)$$
with polynomials
$P_{ik}^{n}(x,v)$, $Q_{ik}^{n}(x,v)$, $R_{ik}^{n}(x,v)$ satisfying
\begin{align}\label{constraint}
	\begin{pmatrix}
		f_{ik}^n\\ v_kf_{ik}^n \\ v_k^2f_{ik}^n
	\end{pmatrix}=\frac{1}{\Delta x \Delta v_i^{n}}\int_{I_{ik}^n} \begin{pmatrix}
	P_{ik}^{n}(x,v)\\Q_{ik}^{n}(x,v)\\R_{ik}^{n}(x,v)
\end{pmatrix} \,dx\,dv, \quad v_k \in \mathcal{G}_i^n
\end{align}
Then, as in \eqref{com rho} and \eqref{com rhouE np1}, we compute
\begin{align*} 
	\begin{split}
		\begin{pmatrix}
			\rho_i^{n+1}\\ \rho_i^{n+1}U_i^{n+1}\\ 2E_i^{n+1}
		\end{pmatrix}&:=\sum_{j\in \mathcal{J}_i^{n+1}} \frac{1}{\Delta x}\int_{\tilde{I}_{i,j}^{n+1}}\begin{pmatrix}
			P^n(x,v)  \\ Q^n(x,v) \\ R^n(x,v)
		\end{pmatrix}\,dx\,dv,\cr
		d_vRT_i^{n+1}&:=(2E_i^{n+1}- \rho_i^{n+1}U_i^{n+1})/\rho_i^{n+1}.
	\end{split}
\end{align*}
where 
$$\tilde{I}_{ij}^{n+1}:=\{(x,v) \in \mathbb{R}^2 |\,  x_{i-\frac{1}{2}}- v \Delta t\leq x \leq  x_{i+\frac{1}{2}} - v \Delta t,\quad   v_{j}-\Delta v_i^{n+1}/2\leq v \leq v_{j}-\Delta v_i^{n+1}/2
\}$$
where $\tilde{I}_{ij}^{n+1}$ denotes the set of the characteristic feet which come from the cell $I_{ij}^{n+1}$. For better understanding of this region, we refer a red parallelogram in Figure \ref{parallelogram np1}.
Then, we obtain
\begin{align*} 
\begin{split}		
	\sum_i	\begin{pmatrix}
		\rho_i^{n+1}\\ \rho_i^{n+1}U_i^{n+1}\\ 2E_i^{n+1}
	\end{pmatrix}\Delta x&=\sum_i\sum_{j\in \mathcal{J}_i^{n+1}}\int_{\tilde{I}_{i,j}^{n+1}}\begin{pmatrix}
		P^n(x,v)  \\ Q^n(x,v) \\ R^n(x,v)
	\end{pmatrix}\,dx\,dv\cr
	&=\sum_i\sum_{k\in \mathcal{J}_i^{n}}\int_{I_{i,k}^{n}}\begin{pmatrix}
	P^n(x,v)  \\ Q^n(x,v) \\ R^n(x,v)
\end{pmatrix}\,dx\,dv=\sum_{ik} \begin{pmatrix}
f_{ik}^n\\ v_kf_{ik}^n \\ v_k^2f_{ik}^n
\end{pmatrix} \Delta x \Delta v_i^n.
\end{split}
\end{align*}
Note that the second equation holds only when the new phase space $\cup_i \mathcal{G}_i^{n+1} $ includes the old one $\cup_i \mathcal{G}_i^{n}$. In practice, here we may introduce small conservation errors due to the truncation of the velocity domain.

\subsection{Construction of piecewise polynomials}\label{sec con poly}
In this section, we describe how we reconstruct polynomials $P_{ik}^n$, $Q_{ik}^n$, $R_{ik}^n$. In particular, our goal is to construct a polynomial of degree one so that the reconstruction technique introduced in \cite{CBRY1,CBRY2} gives third order accuracy for smooth solutions.
Let us consider polynomials of the following forms:
\begin{align*}
	P_{ik}^{n}(x,v)&=f_{ik}^n + f_{ik}^{n\prime}(x-x_i) + f_{ik}^{n\backprime}(v-v_k)
	\cr
	Q_{ik}^{n}(x,v)&=v_kf_{ik}^n + (vf)_{ik}^{n\prime}(x-x_i) + (vf)_{ik}^{n\backprime}(v-v_k)\cr
	R_{ik}^{n}(x,v)&=v_k^2f_{ik}^n + (v^2f)_{ik}^{n\prime}(x-x_i) + (v^2f)_{ik}^{n\backprime}(v-v_k),
\end{align*}
where $f_{ik}^{n\prime}$ and $f_{ik}^{n\backprime}$ are approximations of the first order derivative with respect to $x$ and $v$ directions. Similarly, other coefficients $(vf)_{ik}^{n\prime}$, $(vf)_{ik}^{n\backprime}$, $(v^2f)_{ik}^{n\prime}$, $(v^2f)_{ik}^{n\backprime}$ are denoted. Then, these three polynomials automatically satisfy \eqref{constraint}.
The first step is to compute the slopes $f_{ik}^{n\backprime}$, $(vf)_{ik}^{n\backprime}$, $(v^2f)_{ik}^{n\backprime}$ with respect to $v$-direction using a modified minmod limiter: 
\begin{align*}
	f_{ik}^{n\backprime}&=\minmod\left(\minmod\left(\theta\frac{f_{i,k}^n-f_{i,k-1}^n}{\Delta v_i^n},\theta\frac{f_{i,k+1}^n-f_{i,k}^n}{\Delta v_i^n} \right), \frac{f_{i,k+1}^n-f_{i,k-1}^n}{2\Delta v_i^n}\right)\cr
	vf_{ik}^{n\backprime}&=\minmod\left(\minmod\left(\theta\frac{v_kf_{i,k}^n-v_{k-1}f_{i,k-1}^n}{\Delta v_i^n},\theta\frac{v_{k+1}f_{i,k+1}^n-v_kf_{i,k}^n}{\Delta v_i^n} \right), \frac{v_{k+1}f_{i,k+1}^n-v_{k-1}f_{i,k-1}^n}{2\Delta v_i^n}\right)\cr
	v^2f_{ik}^{n\backprime}&=\minmod\left(\minmod\left(\theta\frac{v_k^2f_{i,k}^n-v_{k-1}^2f_{i,k-1}^n}{\Delta v_i^n},\theta\frac{v_{k+1}^2f_{i,k+1}^n-v_k^2f_{i,k}^n}{\Delta v_i^n} \right), \frac{v_{k+1}^2f_{i,k+1}^n-v_{k-1}^2f_{i,k-1}^n}{2\Delta v_i^n}\right)
\end{align*}
where
\begin{align*}
	\minmod(a,b)=\begin{cases}
		a, \quad \text{if} \quad|a|\leq |b|\\
		b, \quad \, \text{if} \quad|a|>|b|
	\end{cases}
\end{align*}
For $\theta \in [1,2]$, this limiter prevents oscillations for monotone sequences and for $\theta=1$
it reduces to a first order approximation. For numerical simulations, we fix $\theta=1.5$.

Next, we move on to the approximation of the slopes $f_{ik}^{n\prime}$, $(vf)_{ik}^{n\prime}$, $(v^2f)_{ik}^{n\prime}$ with respect to $x$-direction. Since $v_k$ may not belong to $G_{i-1}^n$ and $G_{i+1}^n$, it is necessary to know the values of $f_{i-1}^n(v_k)$, $(vf)_{i-1}^n(v_k)$ and $(v^2f)_{i-1}^n(v_k)$. We approximate the values using
\begin{align*}
	f_{i-1}^n(v_k)&= \frac{1}{\Delta x \Delta v_{i-1}^n}\int_{I_{i-1,k}}P_{i-1,k}^{n}(x,v) \,dx\,dv = \frac{1}{ \Delta v_{i-1}^n}\int_{v_k-\Delta v_{i-1}^{n}/2}^{v_k+\Delta v_{i-1}^{n}/2} f_{i-1,k}^n + f_{i-1,k}^{n\backprime}(v-v_k) \,dx\,dv\cr
	(vf)_{i-1}^n(v_k)&= \frac{1}{\Delta x \Delta v_{i-1}^n}\int_{I_{i-1,k}}Q_{i-1,k}^{n}(x,v) \,dx\,dv = \frac{1}{ \Delta v_{i-1}^n}\int_{v_k-\Delta v_{i-1}^{n}/2}^{v_k+\Delta v_{i-1}^{n}/2} v_kf_{i-1,k}^n + (vf)_{i-1,k}^{n\backprime}(v-v_k) \,dx\,dv\cr
	(v^2f)_{i-1}^n(v_k)&= \frac{1}{\Delta x \Delta v_{i-1}^n}\int_{I_{i-1,k}}R_{i-1,k}^{n}(x,v) \,dx\,dv = \frac{1}{ \Delta v_{i-1}^n}\int_{v_k-\Delta v_{i-1}^{n}/2}^{v_k+\Delta v_{i-1}^{n}/2} v_k^2f_{i-1,k}^n + (v^2f)_{i-1,k}^{n\backprime}(v-v_k) \,dx\,dv.
\end{align*}
Now, we compute $f_{ik}^{n\prime}$, $(vf)_{ik}^{n\prime}$, $(v^2f)_{ik}^{n\prime}$ with a modified minmod limiter: 
\begin{align*}
	f_{ik}^{n\prime}&=\minmod\left(\minmod\left(\theta\frac{f_{i,k}^n-f_{i-1}^n(v_k)}{\Delta x},\theta\frac{f_{i+1}^n(v_k)-f_{i,k}^n}{\Delta x} \right), \frac{f_{i+1}^n(v_k)-f_{i-1}^n(v_k)}{2\Delta x}\right)\cr
	(vf)_{ik}^{n\prime}&=\minmod\left(\minmod\left(\theta\frac{(vf)_{i,k}^n-(vf)_{i-1}^n(v_k)}{\Delta x},\theta\frac{(vf)_{i+1}^n(v_k)-(vf)_{i,k}^n}{\Delta x} \right), \frac{(vf)_{i+1}^n(v_k)-(vf)_{i-1}^n(v_k)}{2\Delta x}\right)\cr
	(v^2f)_{ik}^{n\prime}&=\minmod\left(\minmod\left(\theta\frac{(v^2f)_{i,k}^n-(v^2f)_{i-1}^n(v_k)}{\Delta x},\theta\frac{f_{i+1}^n(v_k)-f_{i,k}^n}{\Delta x} \right), \frac{(v^2f)_{i+1}^n(v_k)-(v^2f)_{i-1}^n(v_k)}{2\Delta x}\right).
\end{align*}
Thus, we obtain the polynomials $P_{ik}^n$, $Q_{ik}^n$, $R_{ik}^n$.

\section{Second order scheme}\label{sec second SL}
\subsection{Second order semi-Lagrangian method for BGK model}
Here we extend the idea of our first order scheme described in Section \ref{sec first SL} to construct a second order scheme. Among various time integrators, we adopt a second order linear multi-step method, the so called backward difference formula (BDF2) in \cite{HW}, which is stable enough to treat the stiffness arising for small Knudsen numbers. Applying the BDF2 method to the problem $y'=f(y,t)$, the method is represented as follows:
\begin{align*} 
	\begin{array}{l}
		\displaystyle	\text{BDF2}: y^{n+1} =\frac{4}{3} y^{n} - \frac{1}{3}y^{n-1} +  \frac{2}{3} \Delta t f(y^{n+1},t_{n+1}).
	\end{array}
\end{align*} 
Now, we describe our second order scheme as follows.

\noindent{\textbf{Step 1: Prediction of macroscopic quantities at $t^{n+1}$ using $\mathcal{G}_i^n$}.} 
We first need to choose velocity nodes for time $t^{n+1}$. For this, we begin by applying BDF2 to the BGK equation in characteristic form, \eqref{SL}:
\begin{align*}
	f^{n+1}(x,v)=\frac{4}{3}f^n(x-v\Delta t,v)-\frac{1}{3}f^{n-1}(x-2v\Delta t,v) +\frac{2}{3} \frac{\Delta t}{\varepsilon} \left(\mathcal{M}^{n+1}(x,v)-f^{n+1}(x,v)\right).
\end{align*}
Based on this,
we precompute $\rho_i^{*}$, $u_i^{*}$ and $E_i^{*}$ and temperature $T_i^{*}$:
\begin{align}\label{com rhouE2}
	\begin{split}
		\begin{pmatrix}
			\rho_i^*\\ \rho_i^*U_i^*\\ 2E_i^*
		\end{pmatrix}&:=\sum_{k\in \mathcal{J}_i^n}\frac{4}{3}\int_{\tilde{I}_{ik}^{n,1}}\begin{pmatrix}
			P^n(x,v)  \\ Q^n(x,v) \\ R^n(x,v)
		\end{pmatrix}\,dx\,dv -\sum_{k\in \mathcal{J}_i^{n-1}}\frac{1}{3}\int_{\tilde{I}_{ik}^{n,2}}\begin{pmatrix}
			P^{n-1}(x,v)  \\ Q^{n-1}(x,v) \\ R^{n-1}(x,v)
		\end{pmatrix}\,dx\,dv,\cr	
		d_vRT_i^*&:=(2E_i^*- \rho_i^*|U_i^*|^2)/\rho_i^*.
	\end{split}
\end{align}
$$\tilde{I}_{ik}^{n,1}:=\{(x,v) \in \mathbb{R}^2 |\,  x_{i-\frac{1}{2}}- v \Delta t\leq x \leq  x_{i+\frac{1}{2}} - v \Delta t,\quad   v_{k}-\Delta v_i^n/2\leq v \leq v_{k}+\Delta v_i^n/2.
\}.$$

$$\tilde{I}_{ik}^{n,2}:=\{(x,v) \in \mathbb{R}^2 |\,  x_{i-\frac{1}{2}}- 2v \Delta t\leq x \leq  x_{i+\frac{1}{2}} - 2v \Delta t,\quad   v_{k}-\Delta v_i^{n-1}/2\leq v \leq v_{k}+\Delta v_i^{n-1}/2.
\}.$$

For high order in space and velocity domain, we construct polynomials $P^n(x,v)$, $Q^n(x,v)$, $R^n(x,v)$ as explained in Section \ref{sec con poly},
while $P^{n-1}(x,v)$, $Q^{n-1}(x,v)$ and $R^{n-1}(x,v)$ are known from the previous step. As we did in the first order scheme, the integration in \eqref{com rhouE2} will be obtained by summation of the integrals for trapezoidal patches using the explicit formula illustrated in \eqref{trap exact integral}. 

\noindent{\textbf{Step 2: Choice of local velocity grids $\mathcal{G}_i^{n+1}$.}
	
	The new grid is defined as in {\textbf{Step 2}} in Section \ref{sec first SL}. This time, however, the $\delta$ is determined by $\delta=\lceil{2\text{CFL}}\rceil+1$. 

\begin{remark}
In case time scales change significantly over time, one would like to use a time step that changes in time as well. 
In this case a a variable stepsize BDF method could be adopted \cite{HWN}.
\end{remark}

%

	\noindent{\textbf{Step 3: Correction of macroscopic quantities at $t^{n+1}$ and computation of $\tilde{f}_{i,j}^{n,1}$, $\tilde{f}_{i,j}^{n,2}$ and $\mathcal{M}_{ij}^{n+1}$.}
	This step improves the prediction of the moments at time $t^{n+1}$. We first compute $\tilde{f}_{i,j}^{n,1}$ and $\tilde{f}_{i,j}^{n,2}$ as approximation of $f(x_i - v_j \Delta t,v_j,t^n)$ and $f(x_i - 2v_j \Delta t,v_j,t^{n-1})$, respectively:
	\begin{align*} 
		\begin{split}
			\tilde{f}_{i,j}^{n,1}&=\frac{1}{\Delta v_i^{n+1}}\int_{\tilde{I}_{ij}^{n+1,1}}P^n(x,v) \,dx\,dv,\cr
			\tilde{f}_{i,j}^{n,2}&=\frac{1}{\Delta v_i^{n+1}}\int_{\tilde{I}_{ij}^{n+1,2}}P^{n-1}(x,v) \,dx\,dv,
		\end{split}
	\end{align*}
	and use this to compute $\rho_i^{n+1}$ (see Figure \ref{BDF2 parallelogram np1}.)
	\begin{align*} 
		\begin{split}
			\rho_i^{n+1}&:=\sum_{j\in \mathcal{J}_i^{n+1}}\left(\frac{4}{3}\tilde{f}_{i,j}^{n,1} - \frac{1}{3}\tilde{f}_{i,j}^{n,2}\right) \Delta v_i^{n+1}.
		\end{split}
	\end{align*}
	The other quantities $U_i^{n+1}$, $E_i^{n+1}$ and temperature $T_i^{n+1}$ are obtained by
	\begin{align*} 
		\begin{split}
			\begin{pmatrix}
				\rho_i^{n+1}U_i^{n+1}\\ 2E_i^{n+1}
			\end{pmatrix}&:=\sum_{j\in \mathcal{J}_i^{n+1}} \left(\frac{4}{3}\int_{\tilde{I}_{ij}^{n+1,1}}\begin{pmatrix}
			Q^n(x,v) \\ R^n(x,v)
		\end{pmatrix}\,dx\,dv - \frac{1}{3}\int_{\tilde{I}_{ij}^{n+1,2}}\begin{pmatrix}
		Q^{n-1}(x,v) \\ R^{n-1}(x,v)
	\end{pmatrix}\,dx\,dv\right),\cr
			d_vRT_i^{n+1}&:=(2E_i^{n+1}- \rho_i^{n+1}U_i^{n+1})/\rho_i^{n+1}.
		\end{split}
	\end{align*}
	Here, we use the piecewise polynomials $P^n$, $Q^n$ and $R^n$ known from \textbf{Step 1}, and the previously stored $P^{n-1}$, $Q^{n-1}$ and $R^{n-1}$ to compute macroscopic moments. Note that we compute moments by integrating these polynomials over the projected region of the grey region along the characteristics (see Figure \ref{BDF2 parallelogram np1}).
\begin{figure}[t]
	\centering
	\begin{subfigure}[t]{1\linewidth}
		\includegraphics[width=1\linewidth]{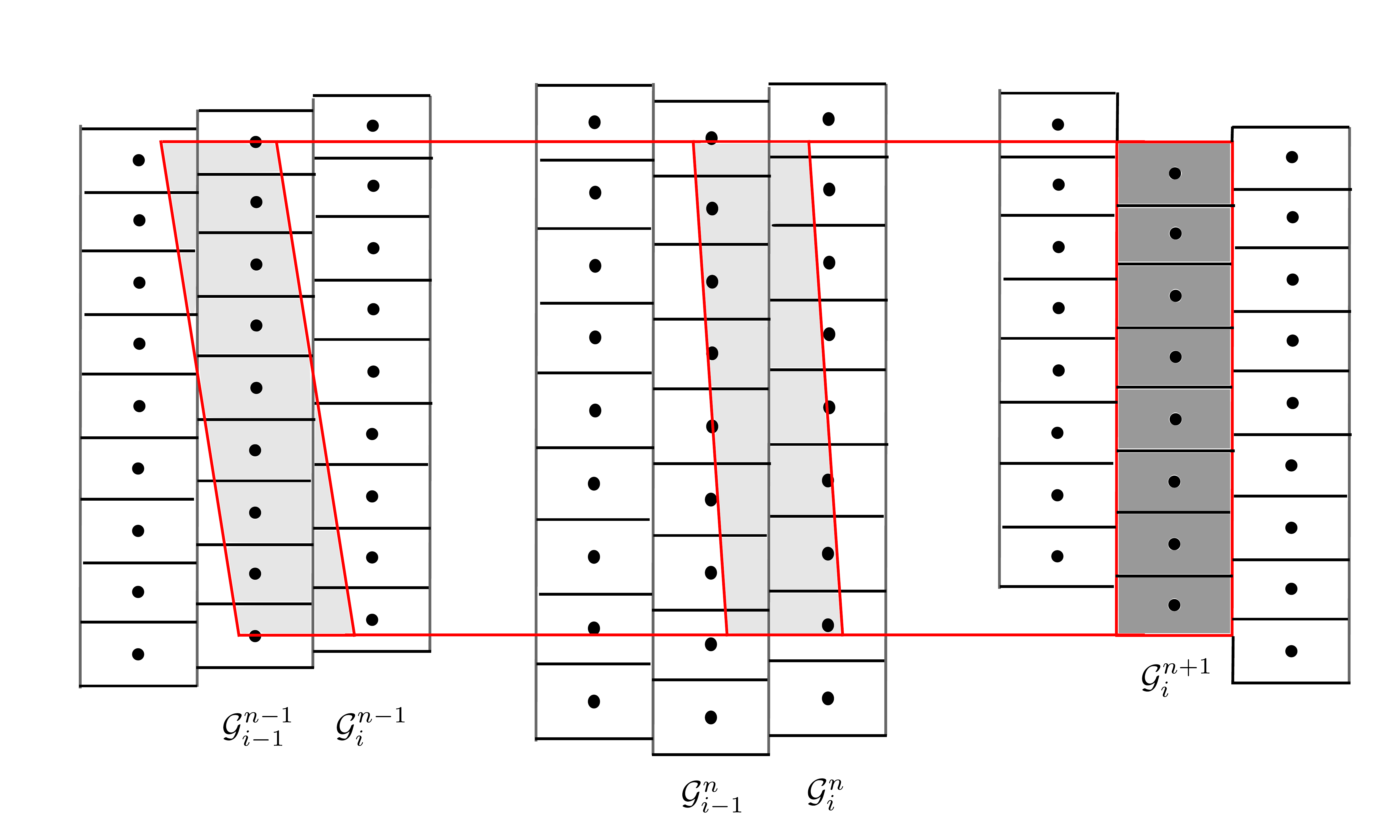}
	\end{subfigure}	
	\caption{
Computation of the 
moments on the new velocity grids for BDF2 (\textbf{Step 3}) obtained by the integration over the light grey regions, which
%
%
%
are the projected domain of the dark grey area along the characteristic curves with velocity belonging to interval  $I_i \times \mathcal{G}_i^{n+1}$. We assume CFL$< 1$ in this figure.}\label{BDF2 parallelogram np1}
	\end{figure}
Next, we compute the local Maxwellian $\mathcal{M}_{ij}^{n+1}$ by
	\begin{align*}
		\mathcal{M}_{ij}^{n+1} =\frac{\rho_{i}^{n+1}}{\sqrt{\left(2 \pi T_{i}^{n+1} \right)^2}}\exp\left(-\frac{|v_{j}-U_{i}^{n+1}|^2}{2T_{i}^{n+1}}\right).
	\end{align*}

	
	\noindent{\textbf{Step 4: Correction using weighted $L^2$-minimization.}
	Now, for each newly defined $v_j\in \mathcal{G}_i^{n+1}$ we attempt to update solutions as
	\begin{align*} 
		f_{i,j}^{n+1}=\tilde{f}_{i,j} + \frac{2}{3}\frac{\Delta t}{\varepsilon} \left(\mathcal{M}_{ij}^{n+1}-f_{i,j}^{n+1}\right).
	\end{align*}
where $\tilde{f}_{i,j}=\frac{4}{3}\tilde{f}_{i,j}^{n,1} - \frac{1}{3}\tilde{f}_{i,j}^{n,2}$. As in Section \ref{sec first SL}, we need to modify $\tilde{f}_{i,j}$ and $\mathcal{M}_{i,j}^{n+1}$ to enforce conservation of discrete moments. Here we also use the weighted $L^2$-minimization technique described in Section \ref{sec weight}, which results in the replacement of each term as follows: 
	\begin{align*}
		\{\mathcal{M}_{i}^{n+1}\}_{j\in \mathcal{J}_i^{n+1}} &\rightarrow M_i^{n+1}:=p_i^{n+1} \circ \mathcal{M}_{i}^{n+1}\cr
		\tilde{f}_{i} &\rightarrow g_{i}:=q_{i}^{n+1} \circ \tilde{f}_{i}.
	\end{align*}
	
	\noindent{\textbf{Step 5: Update of solution $f_{ij}^{n+1}$}.}
	The final step is to update solution using 
	\begin{align*}
		f_{i,j}^{n+1}=g_{i,j}+ \frac{\Delta t}{\varepsilon} \left(M_{ij}^{n+1}-f_{i,j}^{n+1}\right).
	\end{align*}

\section{Numerical tests}\label{sec: numeircal test}
In this section, we perform several tests
checking the accuracy, efficiency and robustness of our scheme. In all tests, time step is obtained by fixing CFL number.
For simplicity, in each run we use the same value of time step through the whole calculation. The time step $\Delta t$ is determined by 
\[
\displaystyle \text{CFL}= \max_i \max_{ v_j \in \mathcal{G}_i^0} |v_j| \frac{ \Delta t}{\Delta x}.
\]

For updating velocity grids in each time step, we use $\alpha=10$, $\beta=0.5$ for all numerical tests.

\subsection{Accuracy test}
	To check the accuracy, we consider the same accuracy test in \cite{GRS}. The initial distribution is given by the Maxwellian:
	\begin{align*} 
		f_0(x,v)=\frac{\rho_0}{\sqrt{2 \pi RT_0}}\exp\left(-\frac{|v-u_0(x)|^2}{2RT_0}\right),\quad \rho_0=T_0=1, x\in [-1,1]
	\end{align*}
	with initial velocity $u_0$:
	\[
	u_0(x) = 0.1 \exp\left(-(10x - 1)^2\right) - 2 \exp\left(-(10x + 3)^2\right).
	\]
with gas constant $R=1$. We impose periodic boundary condition on the interval $ x \in [-1, 1]$, and velocity domain $ v \in [-10, 10]$ up to final time $T_f=0.32$. To check accuracy we fix CFL$=2.4$ and take different mesh sizes in space, $N_x=40, 80, 160, 320$. For velocity, we initially set $N_v=60$ grid points, and local velocity approach is applied for time step $n>1$.
In Table \ref{tab acc}, we reduce the size of space and velocity grid together. In Table \ref{tab acc2}, we fix the number of velocity nodes and reduce the size of mesh in space. In Tables \ref{tab acc}-\ref{tab acc2}, we confirm that the proposed scheme attain order 2 for various $\varepsilon$.


%
%

\begin{table}[t]
	\centering
	{\begin{tabular}{|cccccccc|}
			\hline
			\multicolumn{8}{ |c| }{Relative $L^1$ error and order of density} \\ \hline
			\multicolumn{1}{ |c }{}&
			\multicolumn{1}{ |c| }{}& \multicolumn{2}{ c  }{$\varepsilon=10^{-6}$} & \multicolumn{2}{ |c }{$\varepsilon=10^{-4}$}& \multicolumn{2}{ |c| }{$\varepsilon=10^{-2}$} \\ \hline
			\multicolumn{1}{ |c }{$(N_v,2N_v)$}&
			\multicolumn{1}{ |c|  }{$(N_x,2N_x)$} &
			\multicolumn{1}{ c  }{error} &
			\multicolumn{1}{ c|  }{rate} &
			\multicolumn{1}{ |c  }{error} &
			\multicolumn{1}{ c  }{rate} &
			\multicolumn{1}{ |c  }{error} &
			\multicolumn{1}{ c|  }{rate}   \\ 
			\hline
			\hline
			\multicolumn{1}{ |c }{(20,40)}&
			\multicolumn{1}{ |c|  }{$(40,80)$}
			&4.5975e-03      
			&1.48
			&4.5509e-03      
			&1.50
			&2.5388e-03   
			&2.45
			\\
			\multicolumn{1}{ |c }{(40,80)}&
			\multicolumn{1}{ |c|  }{$(80,160)$}
			&1.6431e-03   
			&2.04
			&1.6099e-03   
			&2.09
			&4.6288e-04   
			&2.03       
			\\
			\multicolumn{1}{ |c }{(80,160)}&
			\multicolumn{1}{ |c|  }{$(160,320)$}&3.9889e-04
			&&3.7914e-04
			&
			&1.1324e-04
			&
			\\
			\hline
			\hline
	\end{tabular}}
	\caption{Accuracy test for the 1D BGK equation. We use CFL$=2.4$.}\label{tab acc}

	{\begin{tabular}{|cccccccc|}
		\hline
		\multicolumn{8}{ |c| }{Relative $L^1$ error and order of density} \\ \hline
		\multicolumn{1}{ |c }{}&
		\multicolumn{1}{ |c| }{}& \multicolumn{2}{ c  }{$\varepsilon=10^{-6}$} & \multicolumn{2}{ |c }{$\varepsilon=10^{-4}$}& \multicolumn{2}{ |c| }{$\varepsilon=10^{-2}$} \\ \hline
		\multicolumn{1}{ |c }{$N_v$}&
		\multicolumn{1}{ |c|  }{$(N_x,2N_x)$} &
		\multicolumn{1}{ c  }{error} &
		\multicolumn{1}{ c|  }{rate} &
		\multicolumn{1}{ |c  }{error} &
		\multicolumn{1}{ c  }{rate} &
		\multicolumn{1}{ |c  }{error} &
		\multicolumn{1}{ c|  }{rate}   \\ 
		\hline
		\hline
		\multicolumn{1}{ |c }{}&
		\multicolumn{1}{ |c|  }{$(40,80)$}&4.2876e-03    
		&1.47
		&4.2374e-03   
		&1.49
		&1.7108e-03   
		&2.00
		\\
		\multicolumn{1}{ |c }{80}&
		\multicolumn{1}{ |c|  }{$(80,160)$}&1.5439e-03   
		&1.97   
		&1.5127e-03   
		&2.01
		&4.2942e-04   
		&2.38
		\\
		\multicolumn{1}{ |c }{}&
		\multicolumn{1}{ |c|  }{$(160,320)$}&3.9307e-04
		&&3.7494e-04	
		&&8.2516e-05&
		\\
		\hline
		\hline
\end{tabular}}
\caption{Accuracy test for the 1D BGK equation for fixed $N_v=80$. We use CFL$=2.4$.}\label{tab acc2}

\end{table}

\subsection{Riemann problems}
In this test, we consider a Riemann problem to confirm that our scheme with local velocity grid approach is able to reproduce the result obtained by classical SL schemes based on global velocity grids. In particular, we aim to show that $L^2$-minimization should be involved for local velocity grid approach.
As in \cite{BRULL201422}, we take
the initial data to be a local Maxwellian with macroscopic quantities:
\begin{align*} 
	(\rho_0,u_0,T_0) = \left\{\begin{array}{lr}
		(0.0001,0,0.00480208), & \text{for } x\le 0.3\\
		(0.0000125,0,0.00384167), & \text{for } x>0.3\\
	\end{array}\right\}
\end{align*}
with gas constant $R=208.1$. We impose freeflow boundary condition on the interval $ x \in [0, 0.6]$, and velocity domain $ v \in [-15, 15]$ upto final time $T_f=7.34\times 10^{-2}$. We take uniform spatial nodes with $N_x=300$, initial velocity grids with $N_v=600$ and fix a time step using CFL$=2$. To reproduce the same result in \cite{BRULL201422}, we take $\varepsilon=\tau:=CT^\omega/\rho$ used in \cite{BRULL201422}. With the purpose of considering different scales of $\varepsilon$, we take different values of $C=1.08 \times 10^{-p}$, $p=7,8,9$, $\omega=-0.19$ so that Knudsen numbers varies from $10^{-2}$ to $10^{-4}$. The choice of $p=9$ is the case in the literature \cite{BRULL201422}.

In Figs \ref{test2 const100}-\ref{test2 const1}, we compare reference solutions based on global grids with the solutions obtained by local velocity grid approaches with and without application of weighted $L^2$-minimization (BDF2+MM+LVG and BDF2+MM+LVG-no-$L^2$). As a reference solution, we use the approach in \cite{CBRY2} by the combination of BDF2 time discretization and the piecewise linear reconstruction with modified minmod limiter (BDF2+MM).
 In Fig. \ref{test2 const100}, we observe big differences between the solution without weighted $L^2$-minimization (black line) and the other solutions. The main source of such error is due to the lack of weighted $L^2$-minimization. On the other hand, small differences are observable between our approach with weighted $L^2$-minimization and the reference solution. Although we didn't report the result here, we observed that this difference is getting smaller as we take more grid points, which implies that the error comes from the accuracy of spatial reconstructions.

On the other results in Figs. \ref{test2 const10}-\ref{test2 const1}, the solutions based on local velocity grid approach with weighted $L^2$-minimization show very good agreement with reference solutions as well.

\begin{figure}[htbp]
	\centering
	\begin{subfigure}[t]{0.42\linewidth}
		\includegraphics[width=1\linewidth]{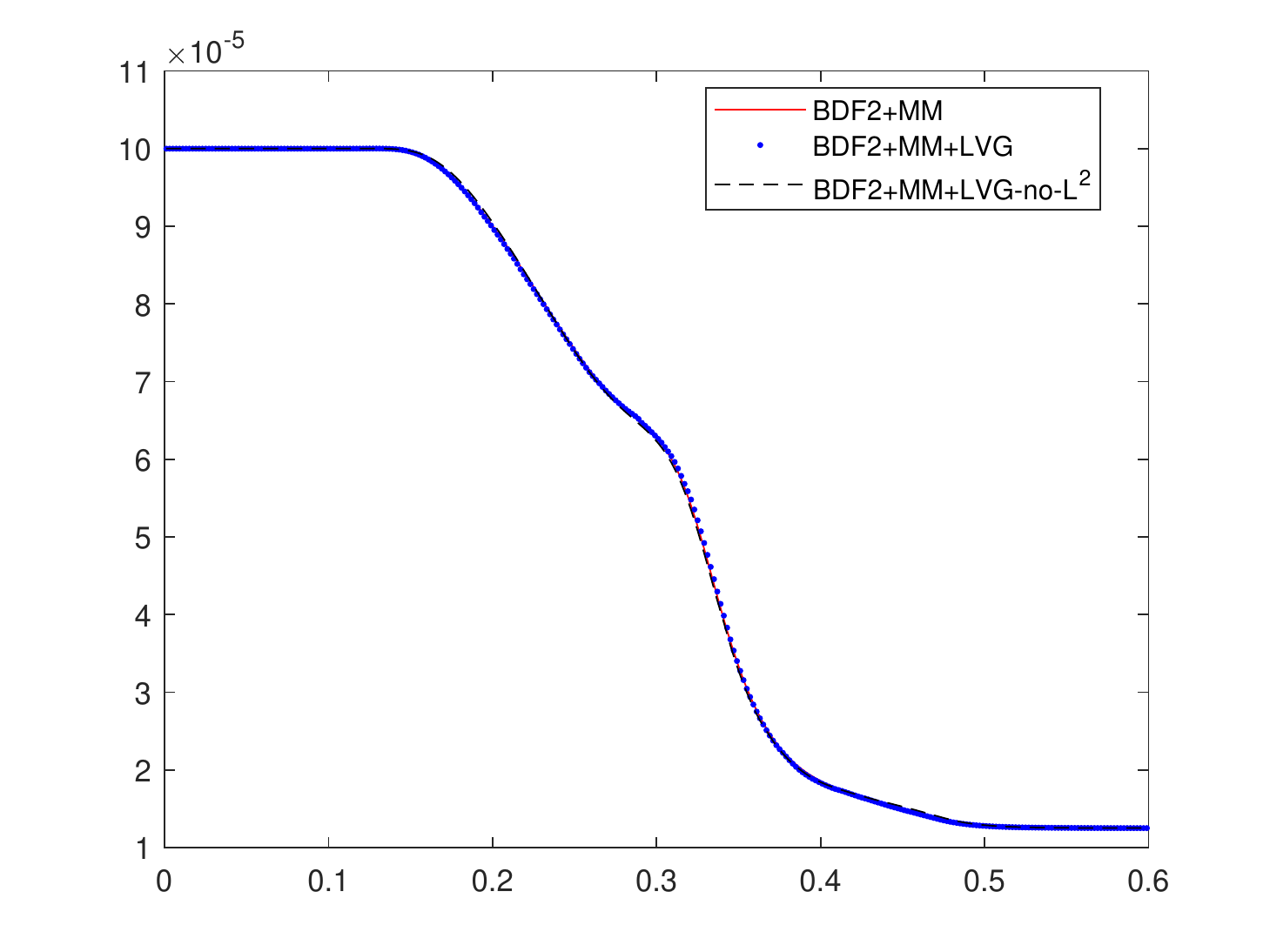}
		\subcaption{Density}
	\end{subfigure}	
	\begin{subfigure}[t]{0.42\linewidth}
		\includegraphics[width=1\linewidth]{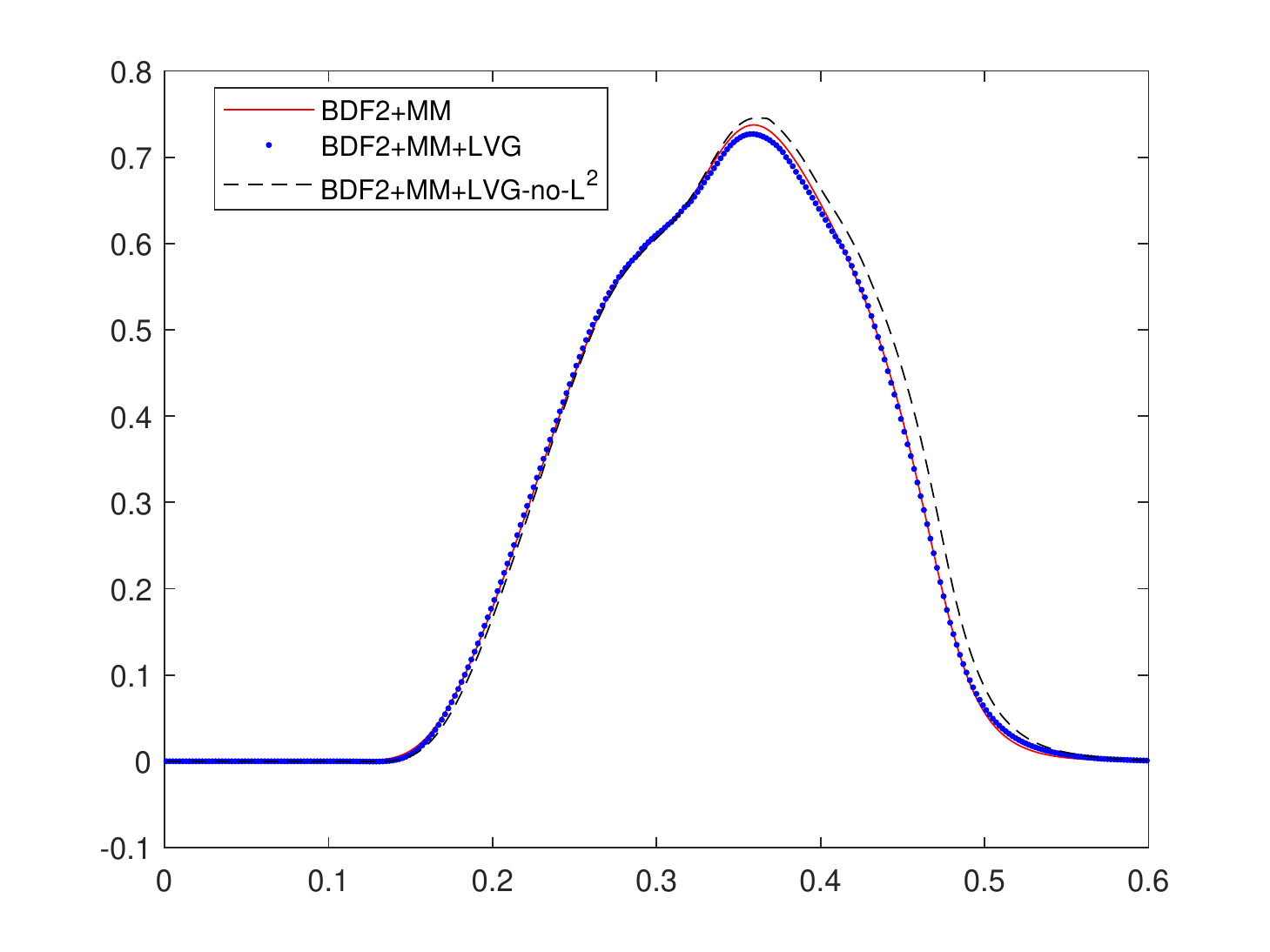}
		\subcaption{Velocity}
	\end{subfigure}	
	\begin{subfigure}[t]{0.42\linewidth}
		\includegraphics[width=1\linewidth]{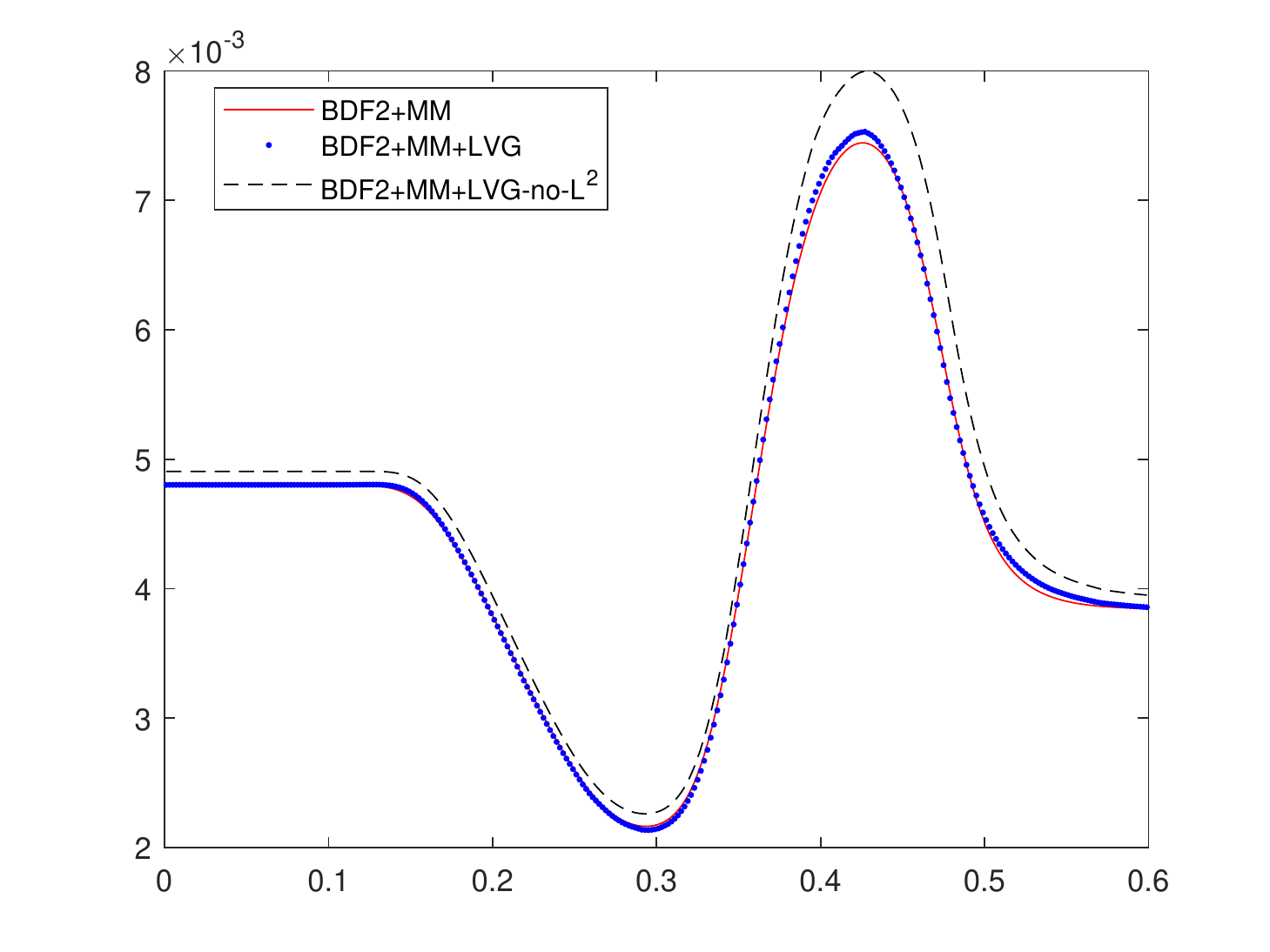}
		\subcaption{Temperature}
	\end{subfigure}	
	\begin{subfigure}[t]{0.42\linewidth}
		\includegraphics[width=1\linewidth]{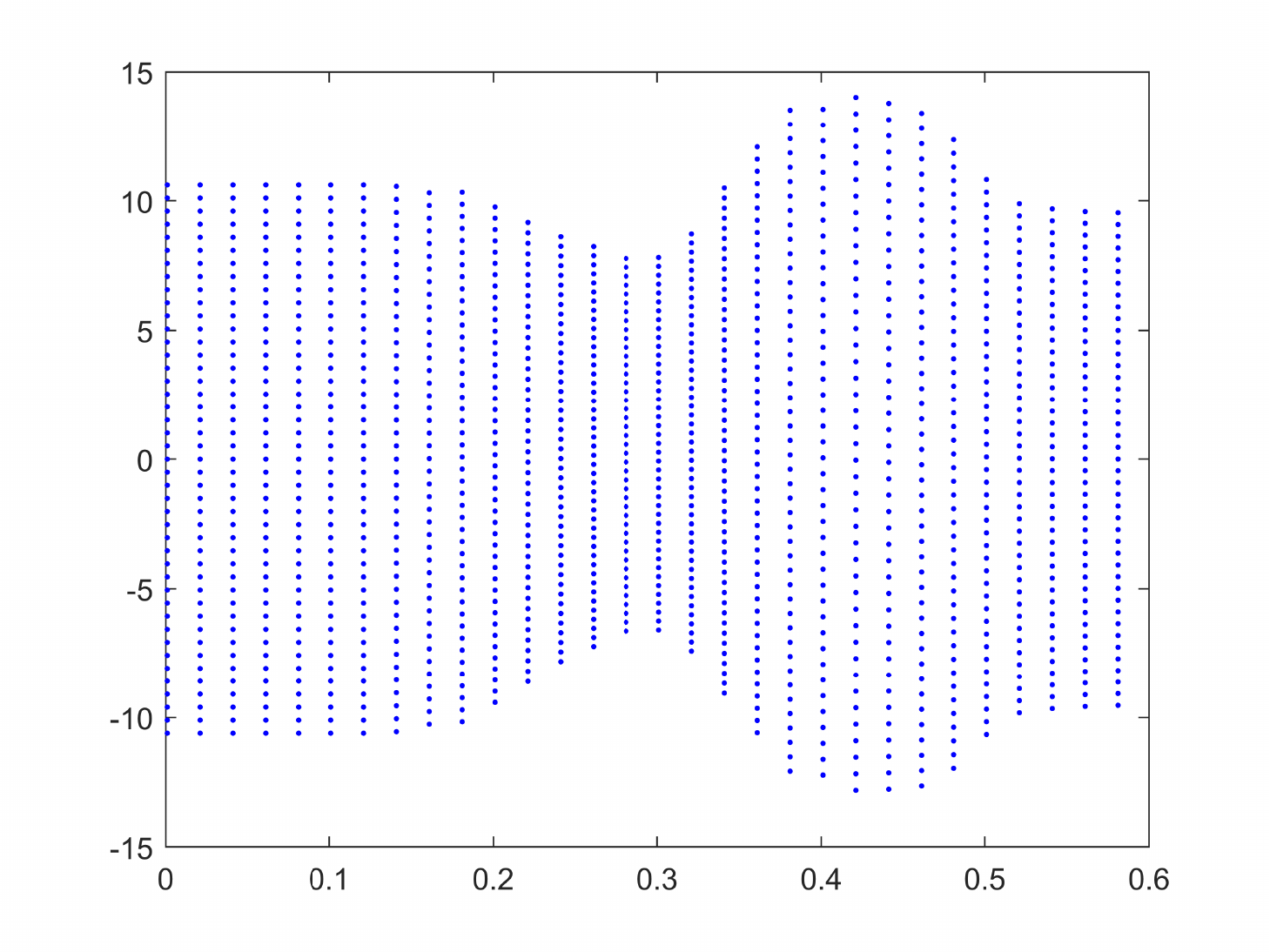}
		\subcaption{Grid points at final time}
	\end{subfigure}	
	
	\begin{subfigure}[t]{0.42\linewidth}
		\includegraphics[width=1\linewidth]{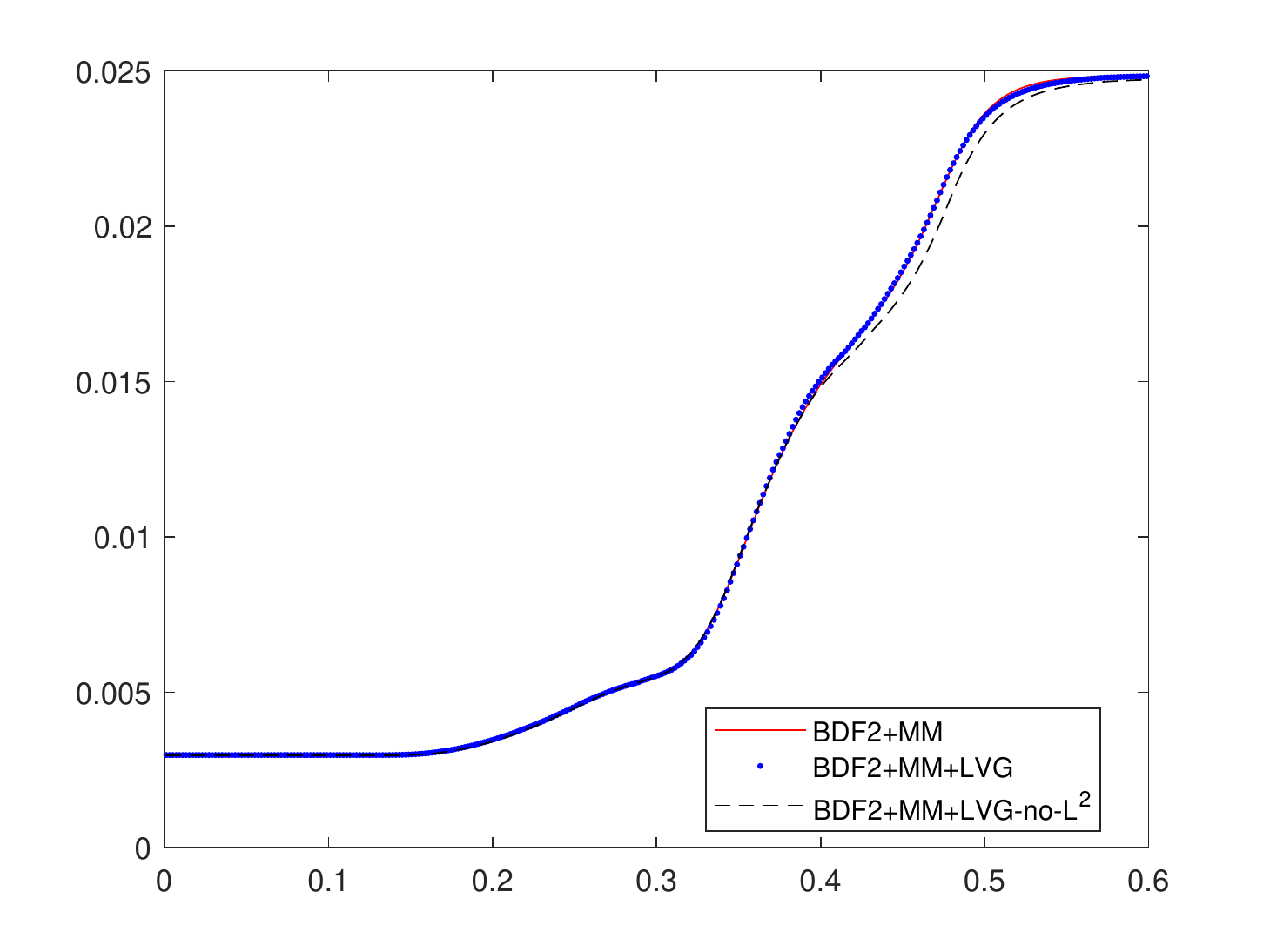}
		\subcaption{$\tau$ at final time}
	\end{subfigure}	
	\caption{Riemann problem with $C=1.08 \times 10^{-7}$. For a classical BDF2+MM scheme (red line) we use $N_v=600$ for each spatial node, while for a local velocity grid approaches we take an average of $N_v=42$.}\label{test2 const100}	
\end{figure}

\begin{figure}[htbp]
	\centering
	\begin{subfigure}[t]{0.42\linewidth}
		\includegraphics[width=1\linewidth]{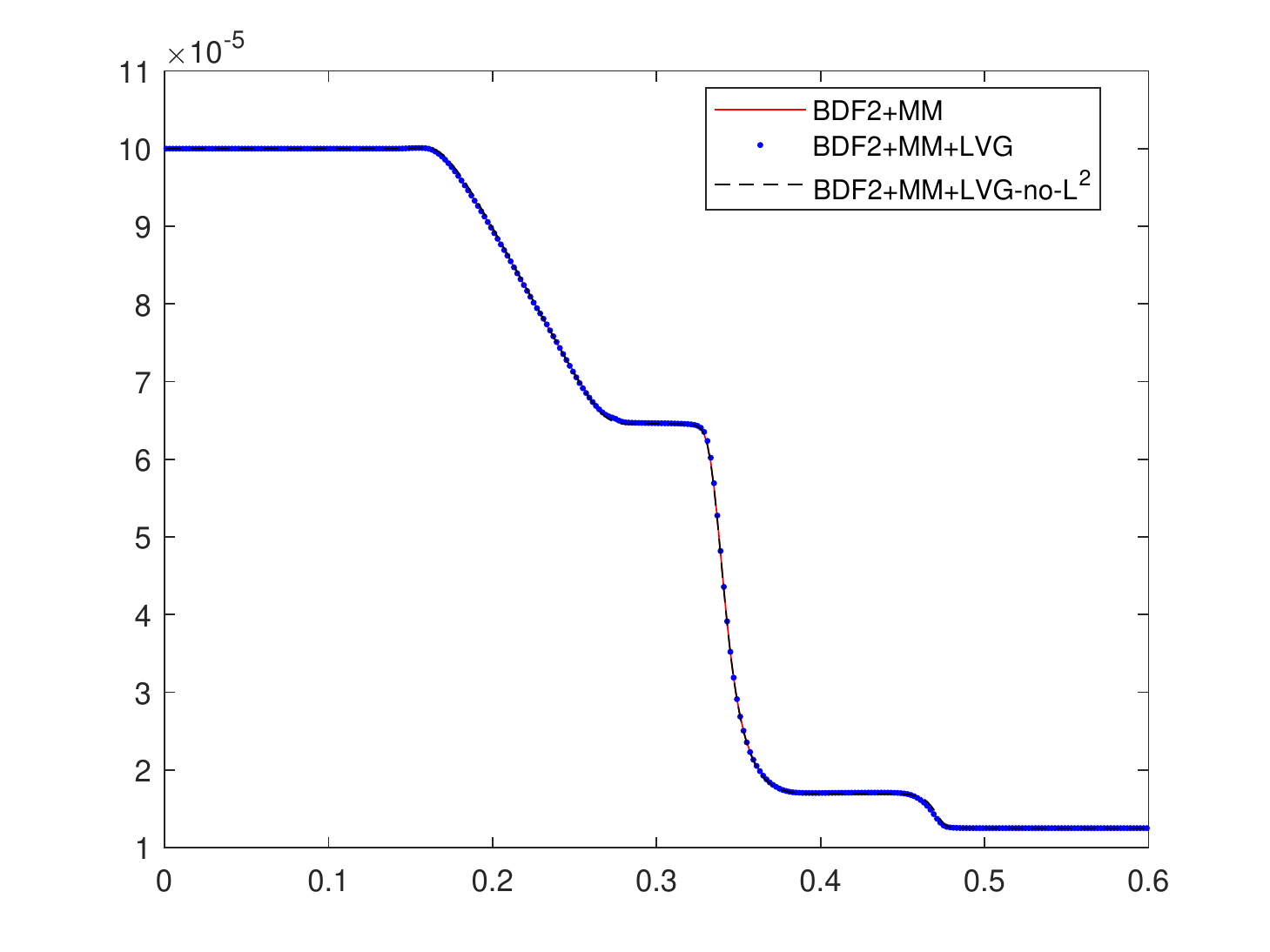}
		\subcaption{Density}
	\end{subfigure}	
	\begin{subfigure}[t]{0.42\linewidth}
		\includegraphics[width=1\linewidth]{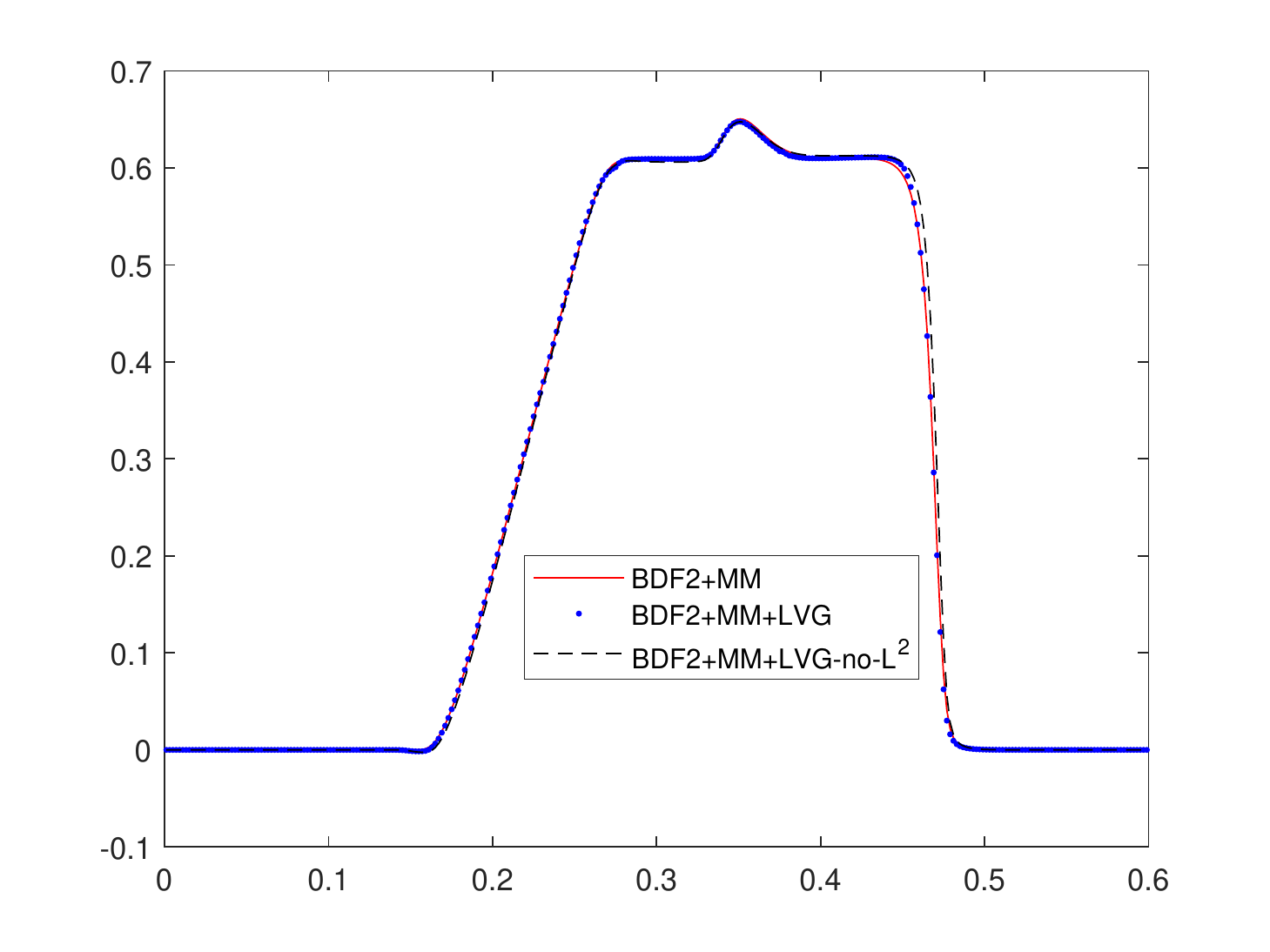}
		\subcaption{Velocity}
	\end{subfigure}	
	\begin{subfigure}[t]{0.42\linewidth}
		\includegraphics[width=1\linewidth]{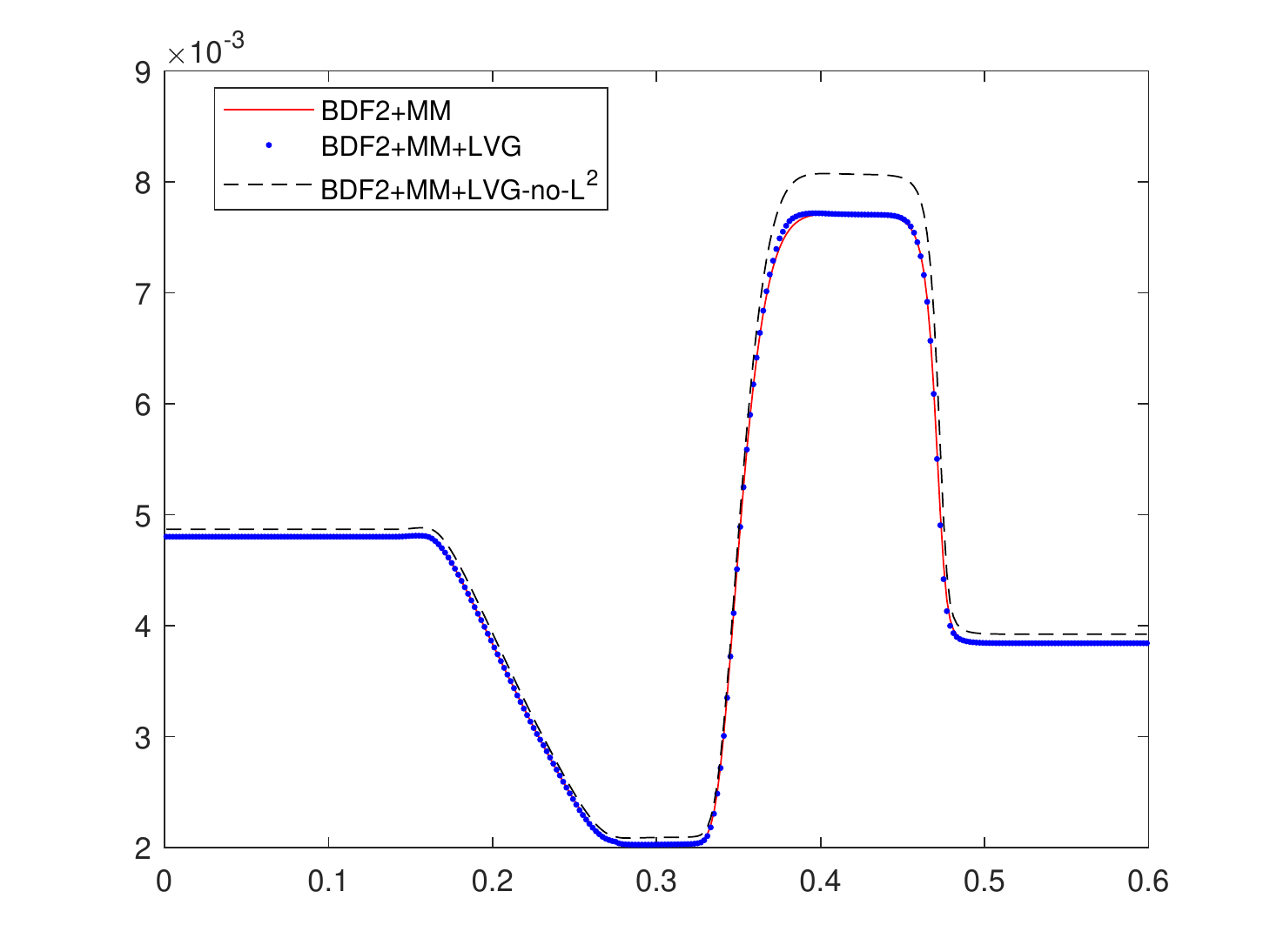}
		\subcaption{Temperature}
	\end{subfigure}	
	\begin{subfigure}[t]{0.42\linewidth}
		\includegraphics[width=1\linewidth]{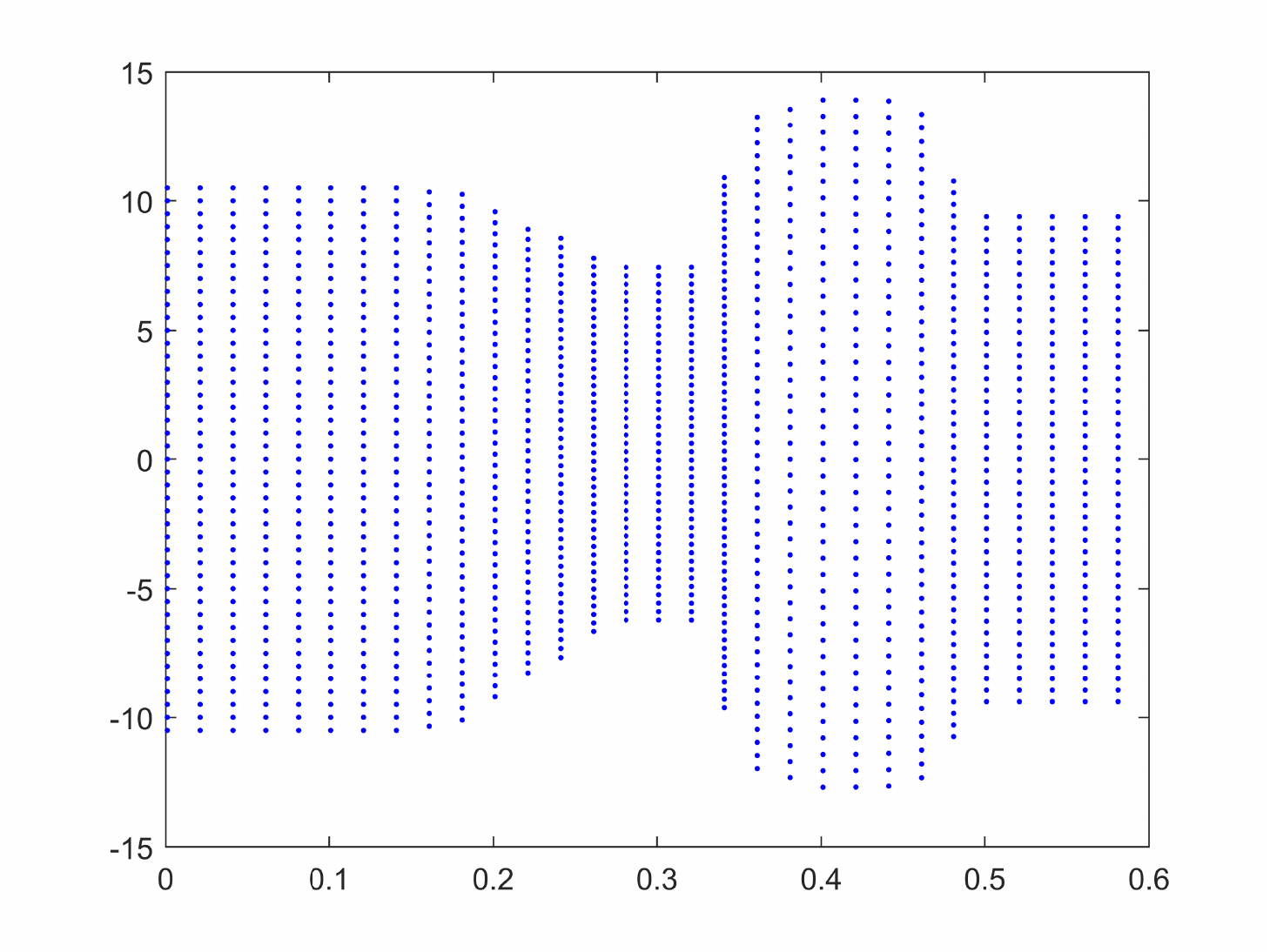}
		\subcaption{Grid points at final time}
	\end{subfigure}	
	
	\begin{subfigure}[t]{0.42\linewidth}
		\includegraphics[width=1\linewidth]{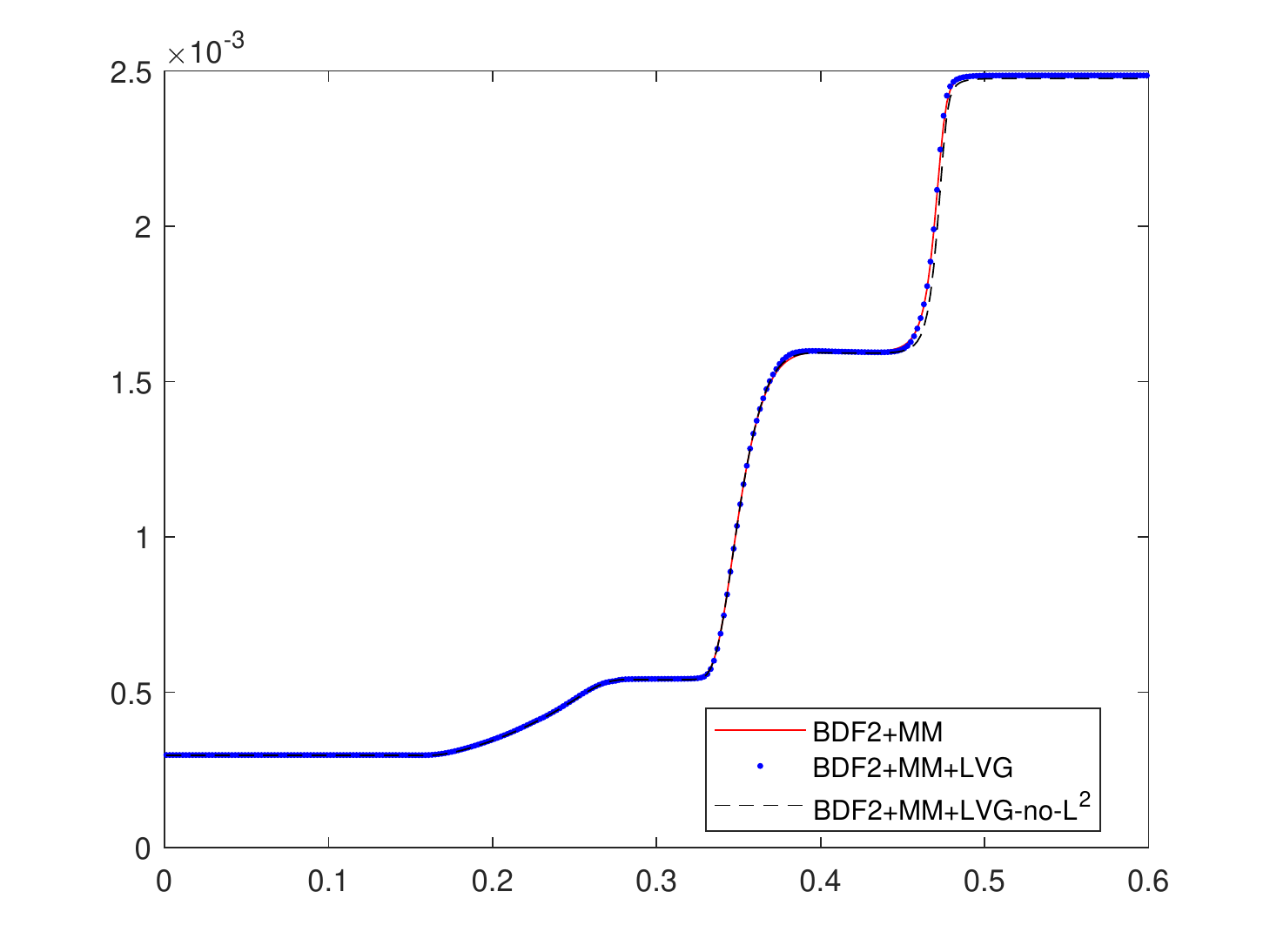}
		\subcaption{$\tau$ at final time}
	\end{subfigure}	
	\caption{Riemann problem with $C=1.08 \times 10^{-8}$ (same as Fig. \ref{test2 const100}).}\label{test2 const10}
	
\end{figure}

\begin{figure}[htbp]
	\centering
	\begin{subfigure}[t]{0.42\linewidth}
		\includegraphics[width=1\linewidth]{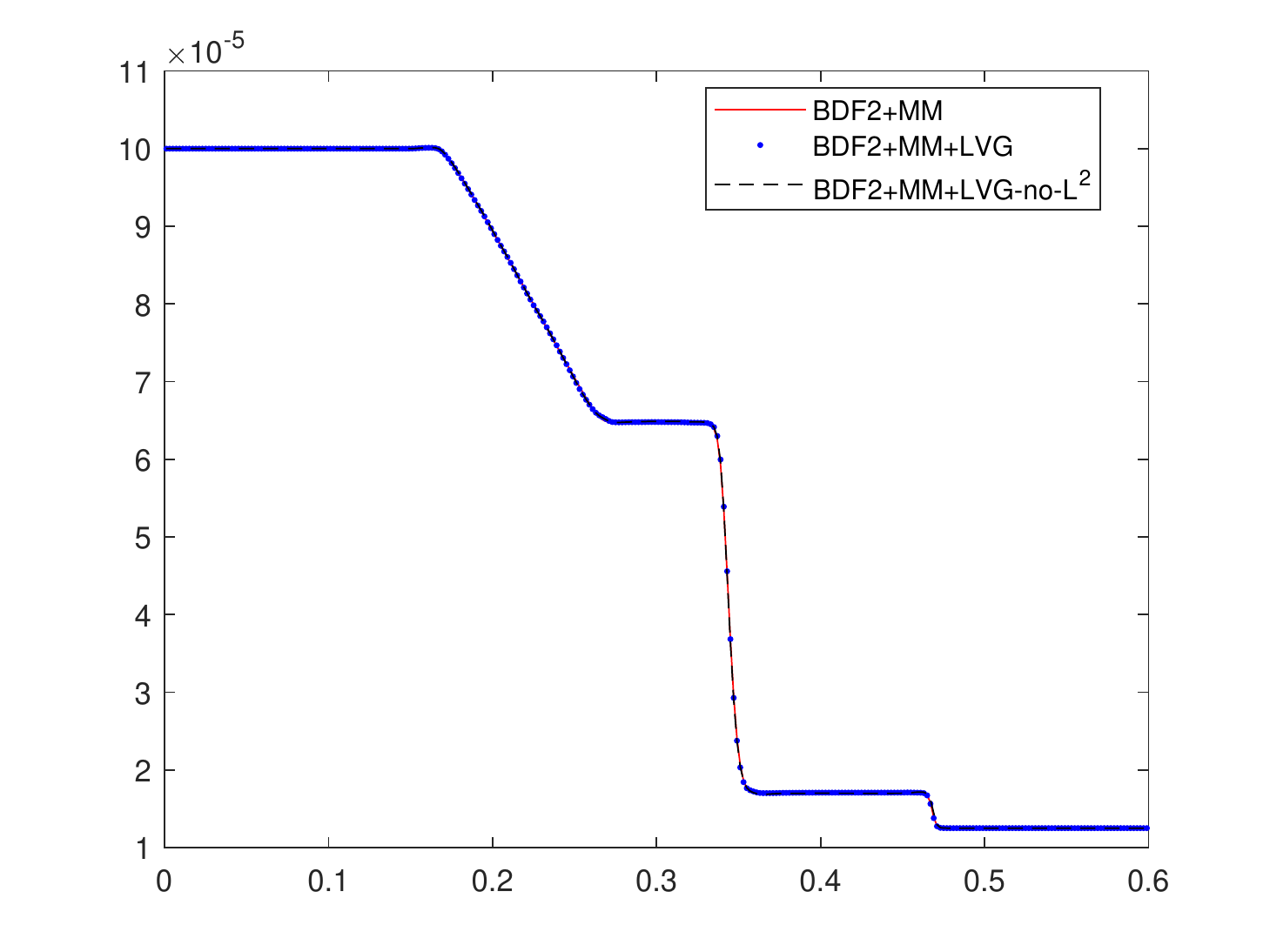}
		\subcaption{Density}
	\end{subfigure}	
	\begin{subfigure}[t]{0.42\linewidth}
		\includegraphics[width=1\linewidth]{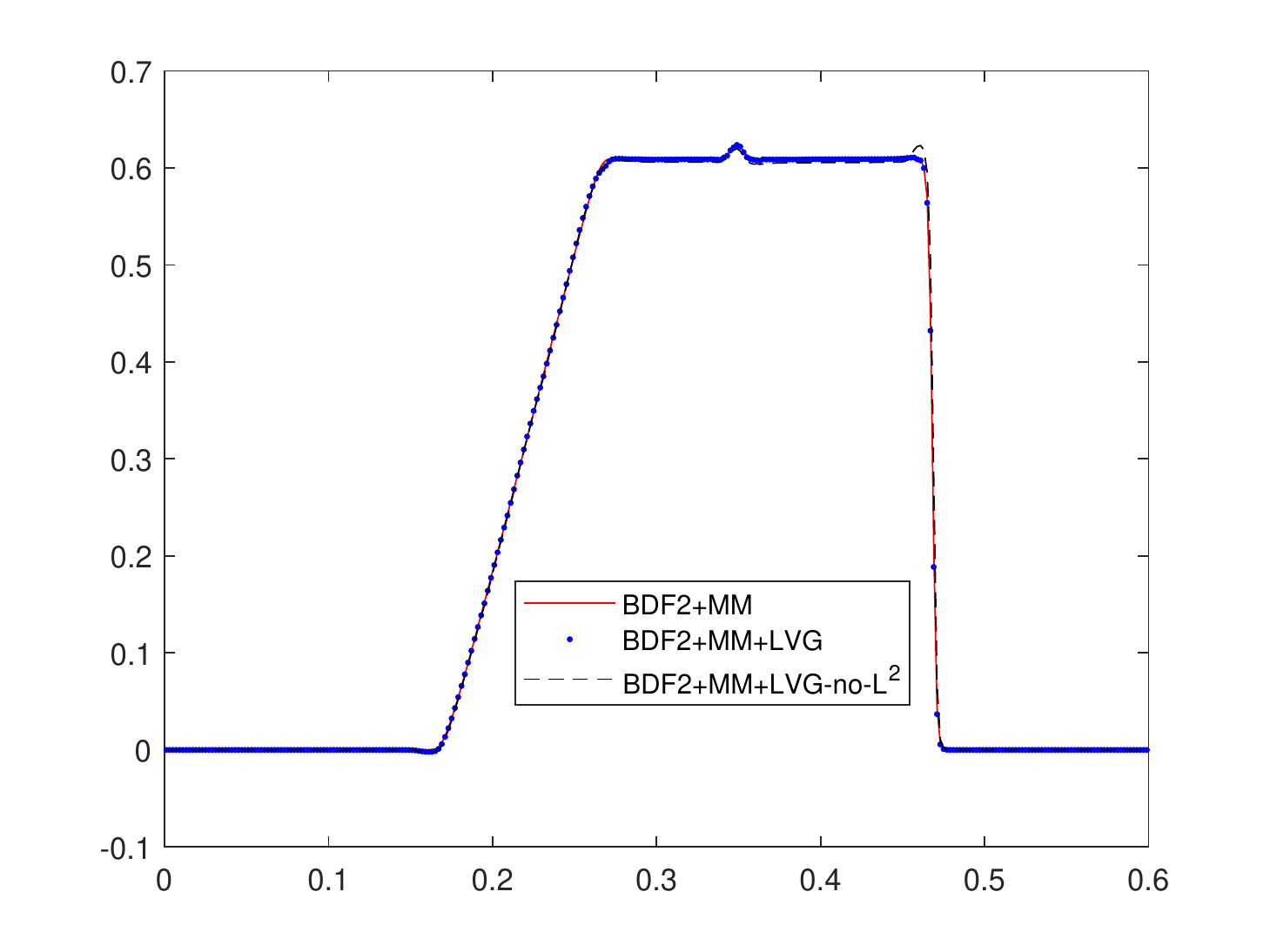}
		\subcaption{Velocity}
	\end{subfigure}	
	\begin{subfigure}[t]{0.42\linewidth}
		\includegraphics[width=1\linewidth]{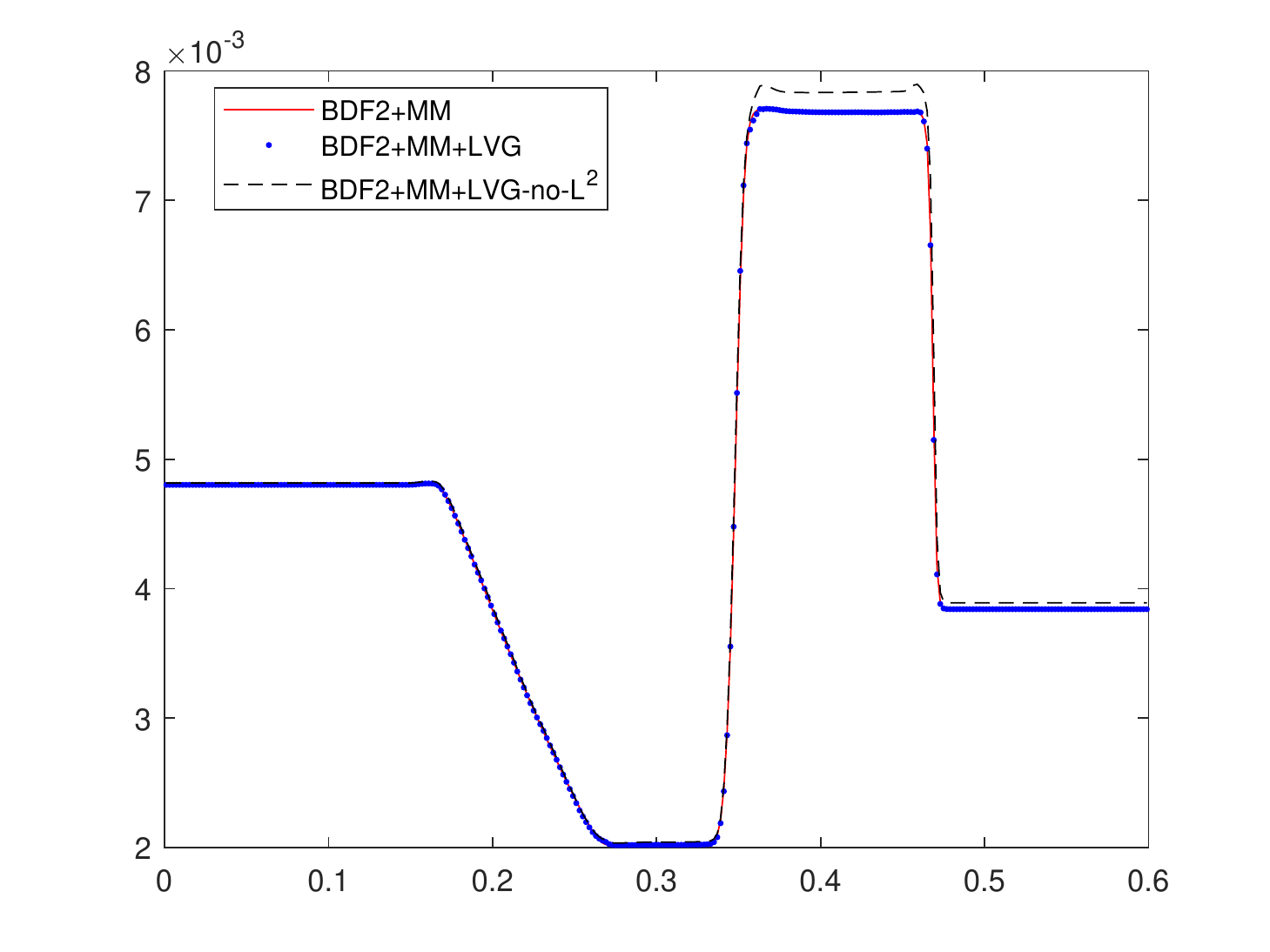}
		\subcaption{Temperature}
	\end{subfigure}	
	\begin{subfigure}[t]{0.42\linewidth}
		\includegraphics[width=1\linewidth]{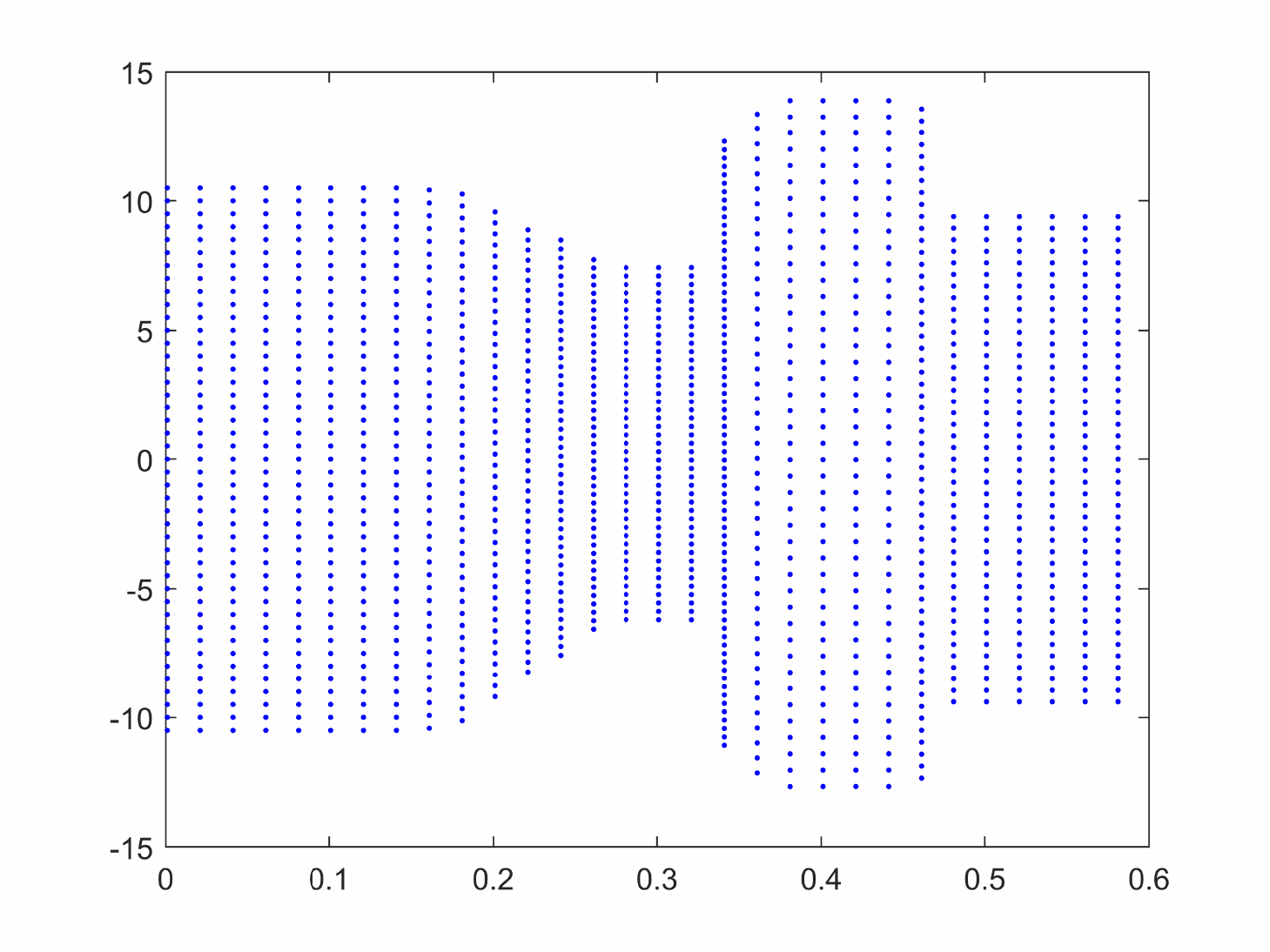}
		\subcaption{Grid points at final time}
	\end{subfigure}	
	
	\begin{subfigure}[t]{0.42\linewidth}
		\includegraphics[width=1\linewidth]{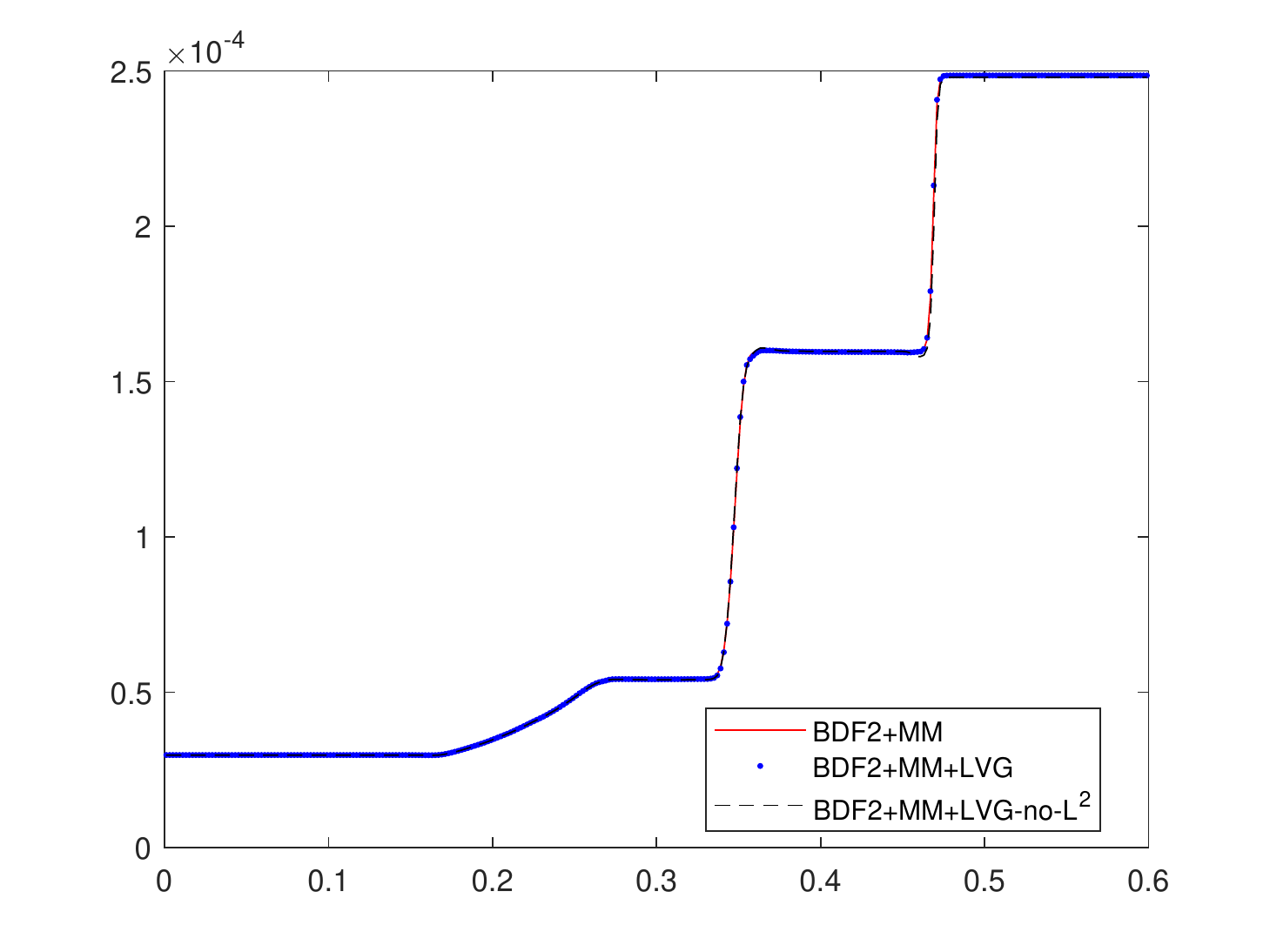}
		\subcaption{$\tau$ at final time}
	\end{subfigure}	
	\caption{Riemann problem with $C=1.08 \times 10^{-9}$ (same as Fig. \ref{test2 const100}).}\label{test2 const1}
\end{figure}

\subsection{Two interacting blast waves}
For a high Mach number, the local velocity approach would be more efficient in that it makes use of fewer grid points. Here we consider a test called ``the two interacting blastwaves"  \cite{BRULL201422}. 
The initial distribution is given by the Maxwellian with initial density $\rho_0(x) \equiv 1 $ and velocity $u_0(x) \equiv 0$. Initial temperature is given by
\begin{align*} 
	T_0(x) = \begin{cases}
		4.8, & \text{for } x\in [0,0.1]\\
		4.8 \times 10^{-5}, & \text{for } x\in [0.1,0.9]\\
		4.8 \times 10^{-1}, & \text{for } x\in [0.9,1]
	\end{cases}
\end{align*}
with gas constant $R=208.1$. We impose freeflow boundary condition on the interval $ x \in [0, 1]$, and compute numerical solutions up to final time $T_f=0.008$. We take uniform grids with $N_x=500$ and fix a time step based on CFL$=2$. To produce the result in  $\varepsilon=\tau:=CT^\omega/\rho$ with $C=1.08 \times 10^{-9}$ as used in \cite{BRULL201422}.

In this problem, the bounds of the global velocity grid of the classical SL schemes are fixed to be $[-190,190]$ to guarantee the conservation up to machine precision. The size of velocity mesh is fixed by $\Delta v=0.1$ which is small enough to satisfy the condition \eqref{dv cond} and enables to resolve Maxwellian distribution corresponding to the smallest initial temperature.

In Fig. \ref{test 3}, we observe that the solution obtained by the local velocity grid approach is able to capture the correct shock position. Note that the use of coarse velocity grids produce some errors, which becomes negligible as the velocity grid is refined. 
Finally, in Fig.\ref{test 3} we depict the local velocity grids used for different space positions.
In this figure, along the vertical direction for each spatial cell, the set of local velocity grids are dotted. Note that the size of mesh in velocity direction varies depending on the temperature. 

\begin{figure}[htbp]
	\centering
	\begin{subfigure}[t]{0.42\linewidth}
		\includegraphics[width=1\linewidth]{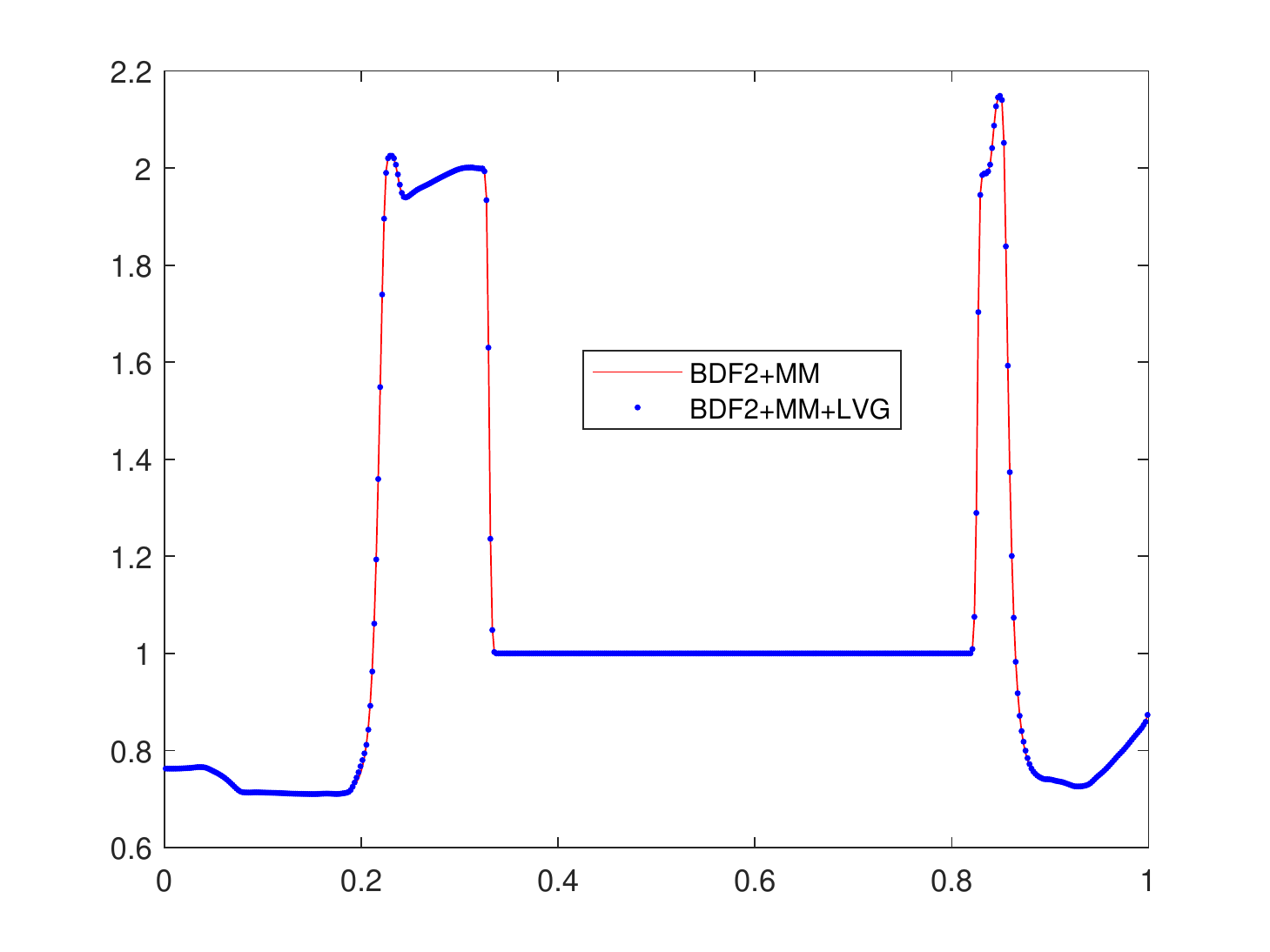}
		\subcaption{Density}
	\end{subfigure}	
	\begin{subfigure}[t]{0.42\linewidth}
		\includegraphics[width=1\linewidth]{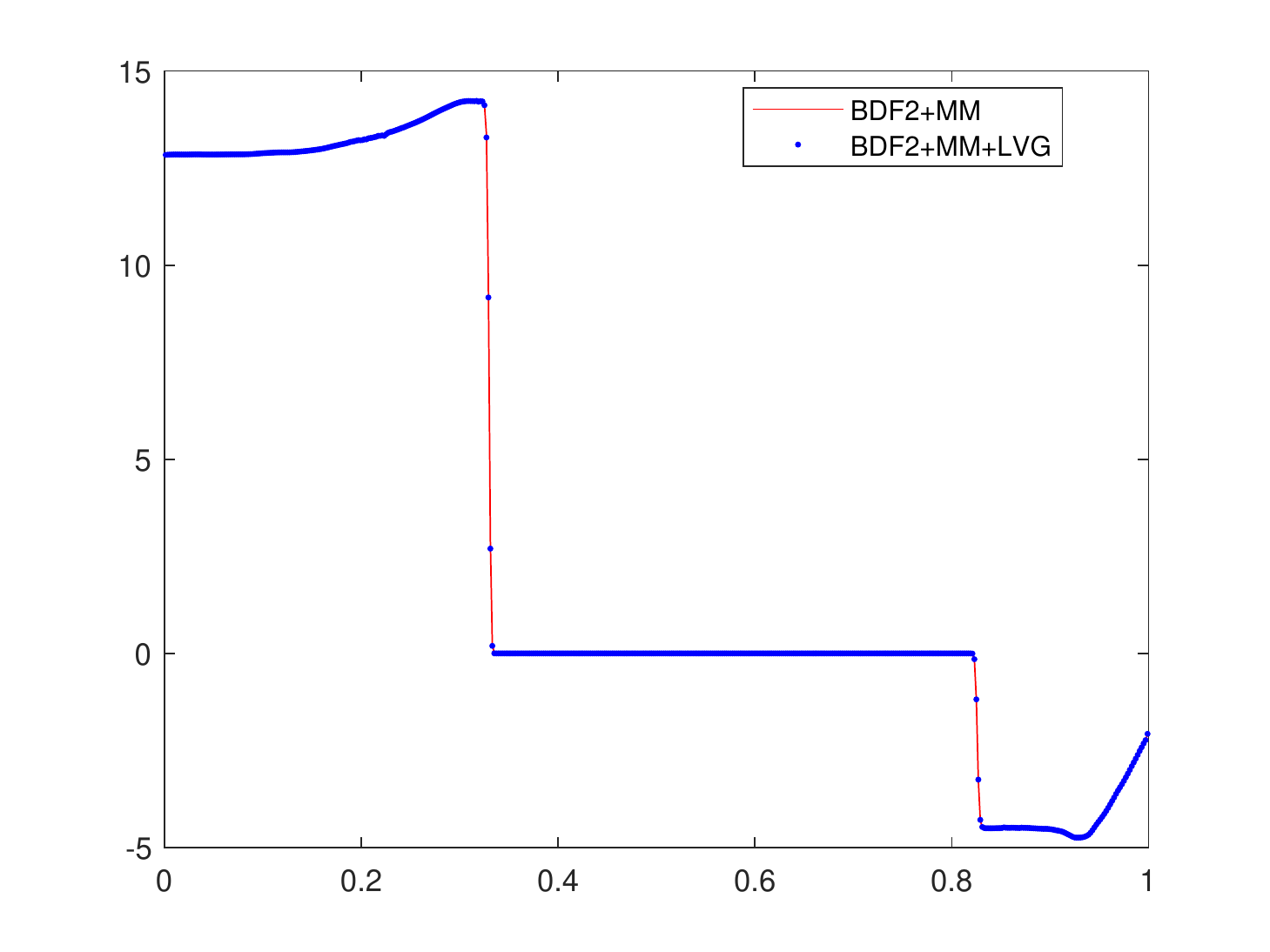}
		\subcaption{Velocity}
	\end{subfigure}	
	\begin{subfigure}[t]{0.42\linewidth}
		\includegraphics[width=1\linewidth]{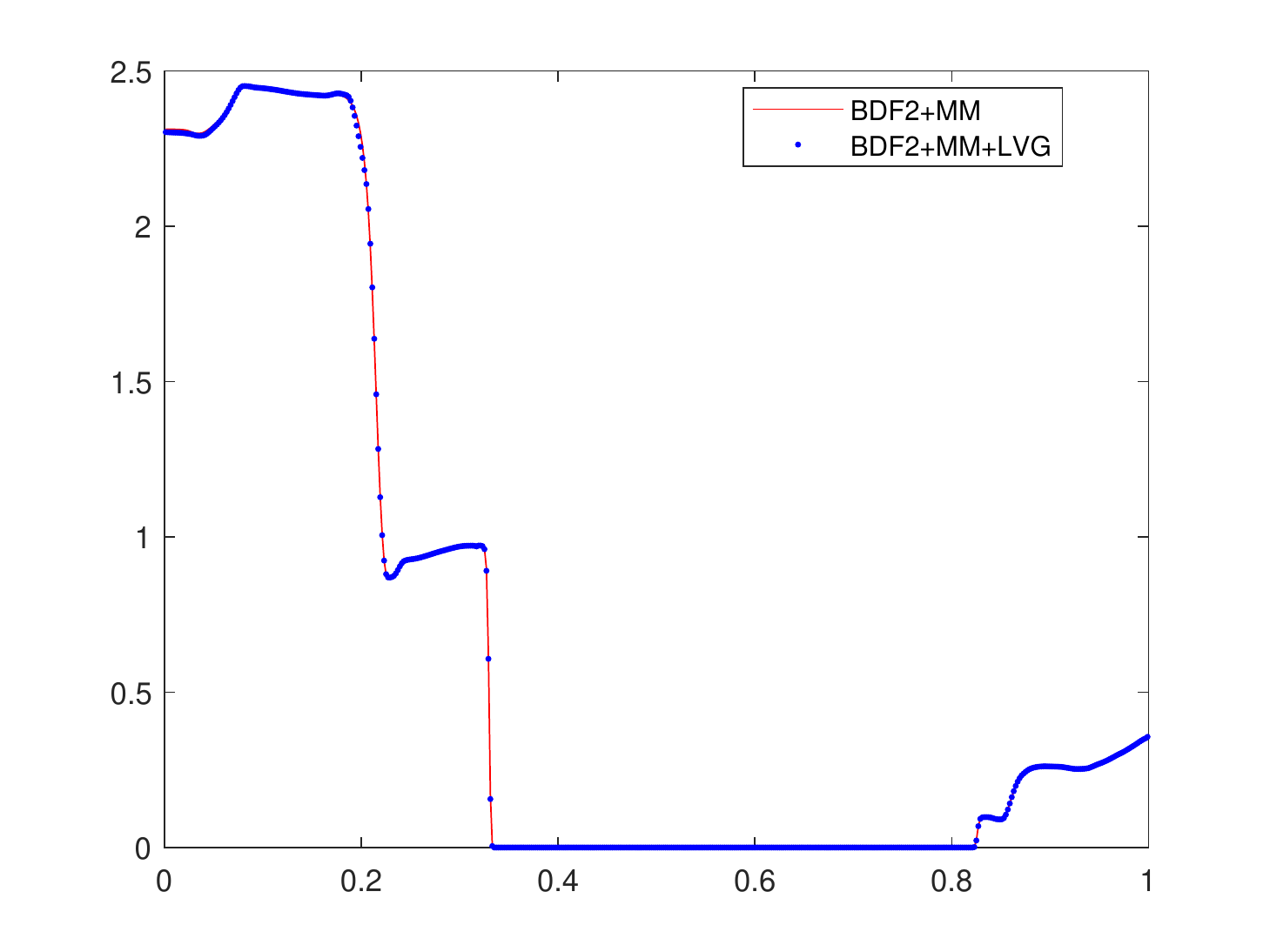}
		\subcaption{Temperature}
	\end{subfigure}	
	\begin{subfigure}[t]{0.42\linewidth}
		\includegraphics[width=1\linewidth]{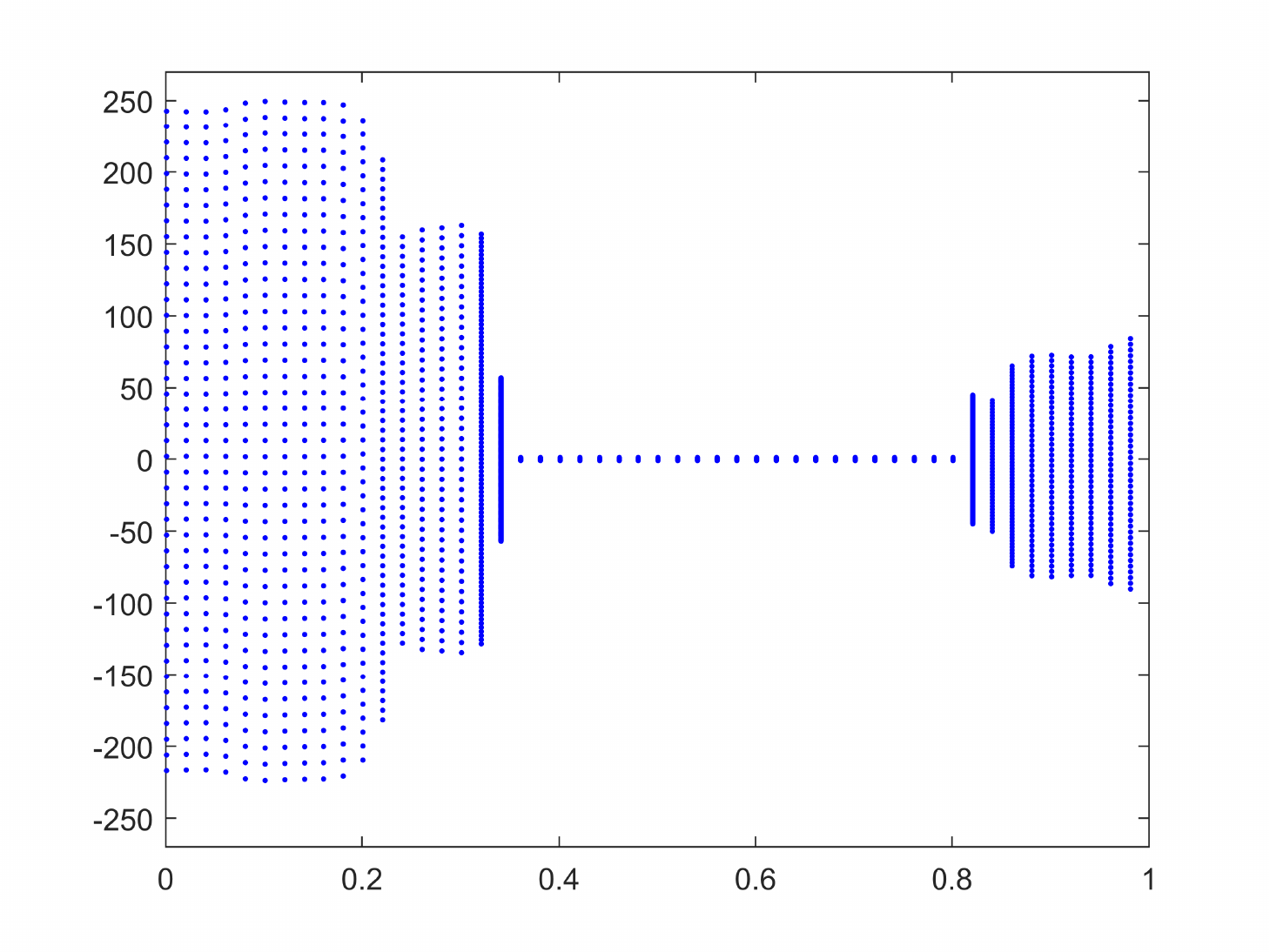}
		\subcaption{Grid points at final time}
	\end{subfigure}	
	\caption{Two interacting blast waves with $C=1.08 \times 10^{-9}$. For a classical BDF2+MM scheme (red line) we use $N_v=3800$ for each spatial node, while for a local velocity grid approaches we take an average of $N_v=42$.}\label{test 3}
\end{figure}

\appendix
\section{Explcit form of integral of polynomials over trapzoidal domain}
In this section, we provide explicit form of the integration of polynomials of degree two over the trapezoidal domain. Let us consider a two dimensional polynomial $p(x,y)$:
\[
p(x,y)= c_{00} + c_{10} x + c_{01}y.
\]
and a trapezoidal domain $\Omega$ 
whose four vertices are given along counterclockwise direction as follows:
\[
(x_1,y_1),\, (x_2,y_2),\, (x_3,y_3),\, (x_4,y_4) \in \mathbb{R}^2.
\]
Denote by $C$ the positively oriented boundary of $\Omega$, and by $S_{ij}$ the line segment of $C$ which connects $(x_i,y_i)$ and $(x_j,y_j)$ (see Figure \ref{parallelogram trap}). Next, we construct a polynomial $q(x,y)$:
\begin{align}\label{q}
	q(x,y)=c_{00}x + \frac{c_{10}}{2} x^2 + c_{01}xy.
\end{align}
which satisfies $\frac{\partial q}{\partial x}=p(x,y)$. Then, the divergence theorem gives
\begin{align*}
	\int_\Omega p(x,y) d xd y&=\int_\Omega \nabla \cdot \begin{pmatrix}
		q(x,y)\\0
	\end{pmatrix}  d xd y \cr
	&=\int_C \begin{pmatrix}
		q(x,y)\\0
	\end{pmatrix} \cdot \textbf{n}(x,y) d s,
\end{align*}
where $\textbf{n}(x,y)$ denotes a outward unit normal vector to the curve $C$. Since the trapezoidal domain is surrounded by four line segments, we split the integral into four parts:
\begin{align*}
	\int_\Omega p(x,y) d xd y
	&=\int_{C_{12}}   \begin{pmatrix}
		q(x,y)\\0
	\end{pmatrix} \cdot \textbf{n}_{12}(x,y) d s + \int_{C_{23}}   \begin{pmatrix}
	q(x,y)\\0
\end{pmatrix}   \cdot \textbf{n}_{23}(x,y) d s\cr
&+\int_{C_{34}}  \begin{pmatrix}
	q(x,y)\\0
\end{pmatrix}  \cdot  \textbf{n}_{34}(x,y) d s + \int_{C_{41}}  \begin{pmatrix}
q(x,y)\\0
\end{pmatrix}  \cdot  \textbf{n}_{41}(x,y) d s.
\end{align*}

\begin{figure}[t]
	\centering
	\begin{subfigure}[t]{0.6\linewidth}
		\includegraphics[width=1\linewidth]{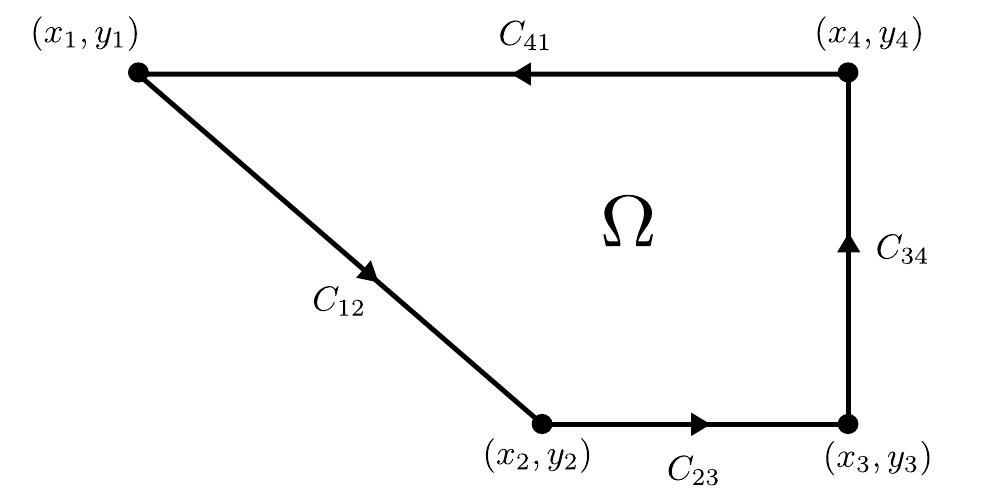}
	\end{subfigure}	
	\caption{Trapzoidal domain $\Omega$ and its positively oriented boundary curve $C= C_{12}\cup C_{23}\cup C_{34}\cup C_{41} $.}\label{parallelogram trap}
\end{figure}

The normal vector $\textbf{n}_{ij}(x,y)$ is given by
\begin{align*}
	\textbf{n}_{ij}(x,y)\equiv  \frac{1}{\sqrt{|x_j-x_i|^2 + |y_j-y_i|^2}} \begin{pmatrix}
		y_j-y_i \\ -(x_j-x_i)
	\end{pmatrix}=: \begin{pmatrix}
	\textbf{n}_{ij}^1\\
	\textbf{n}_{ij}^2
\end{pmatrix}.
\end{align*}
By parametrizing each $C_{ij}$ as 
\begin{align*}
	\textbf{r}_{ij}(t)= \left(x_i + t\Delta x_{ij}, y_i + t\Delta y_{ij}\right),
\end{align*}
with $\Delta x_{ij}=x_j-x_i ,\,\Delta y_{ij}=y_j-y_i$, we obtain
\begin{align*}
	\int_{C_{ij}}   \begin{pmatrix}
		q(x,y)\\0
	\end{pmatrix}   \cdot \textbf{n}_{ij}(x,y) d s&=\int_0^1  q\left(\textbf{r}_{ij}(t)\right) \textbf{n}_{ij}^1(x,y) \left\|\textbf{r}_{ij}'(t)\right\|dt.
\end{align*}
Then, we get
\begin{align*}
	\int_{C_{ij}} 
	\begin{pmatrix}
		q(x,y)\\0
	\end{pmatrix}   \cdot \textbf{n}_{ij}(x,y) d s&=\int_0^1  q\left(\textbf{r}_{ij}(t)\right) \textbf{n}_{ij}^1(x,y) \left\|\textbf{r}_{ij}'(t)\right\|d t\cr
	&=\Delta y_{ij}\int_0^1  q\left(\textbf{r}_{ij}(t)\right)  d t.
\end{align*}
Note that the form of $q(x,y)$ in \eqref{q} gives
\begin{align*}
	\int_0^1  q\left(\textbf{r}_{ij}(t)\right)  d t
	&=
	c_{00}(x_i + \frac{1}{2} \Delta x_{ij})\cr
	&+
	\frac{c_{10}}{2}\left(
	x_i^2 + x_i \Delta x_{ij} +  \frac{1}{3}\Delta x_{ij}^2\right)\cr
	&+
	c_{01}\left(x_iy_i +\frac{x_i\Delta y_{ij} + y_i \Delta x_{ij}}{2}+  \frac{1}{3}\Delta x_{ij} \Delta y_{ij}\right).
\end{align*}
Consequently, we obtain
\begin{align}\label{trap exact integral}
	\begin{split}
		\int_\Omega p(x,y) d xd y
		&=\sum_{(i,j)=(1,2),(2,3),(3,4),(4,1)} \Delta y_{ij}\bigg[
		c_{00}(x_i + \frac{1}{2} \Delta x_{ij})\cr
		&\quad\quad\quad\qquad\qquad\qquad\qquad\qquad+
		\frac{c_{10}}{2}\left(
		x_i^2 + x_i \Delta x_{ij} +  \frac{1}{3}\Delta x_{ij}^2\right)\cr
		&\quad\quad\quad\qquad\qquad\qquad\qquad\qquad+
		c_{01}\left(x_iy_i +\frac{x_i\Delta y_{ij} + y_i \Delta x_{ij}}{2}+  \frac{1}{3}\Delta x_{ij} \Delta y_{ij}\right)\bigg].
	\end{split}
\end{align}
Note that $\Delta y_{23}$ and $\Delta y_{41}$ are zero in our problem. Also, In the case when we integrate polynomials over the grey region shown in Figure \ref{parallelogram intergration}, we consider the integral over $L\cup C \cup R$ and subtract two integral quantities on $L$ and $R$ based on \eqref{trap exact integral}.

\begin{figure}[t]
	\centering
	\begin{subfigure}[t]{0.6\linewidth}
		\includegraphics[width=1\linewidth]{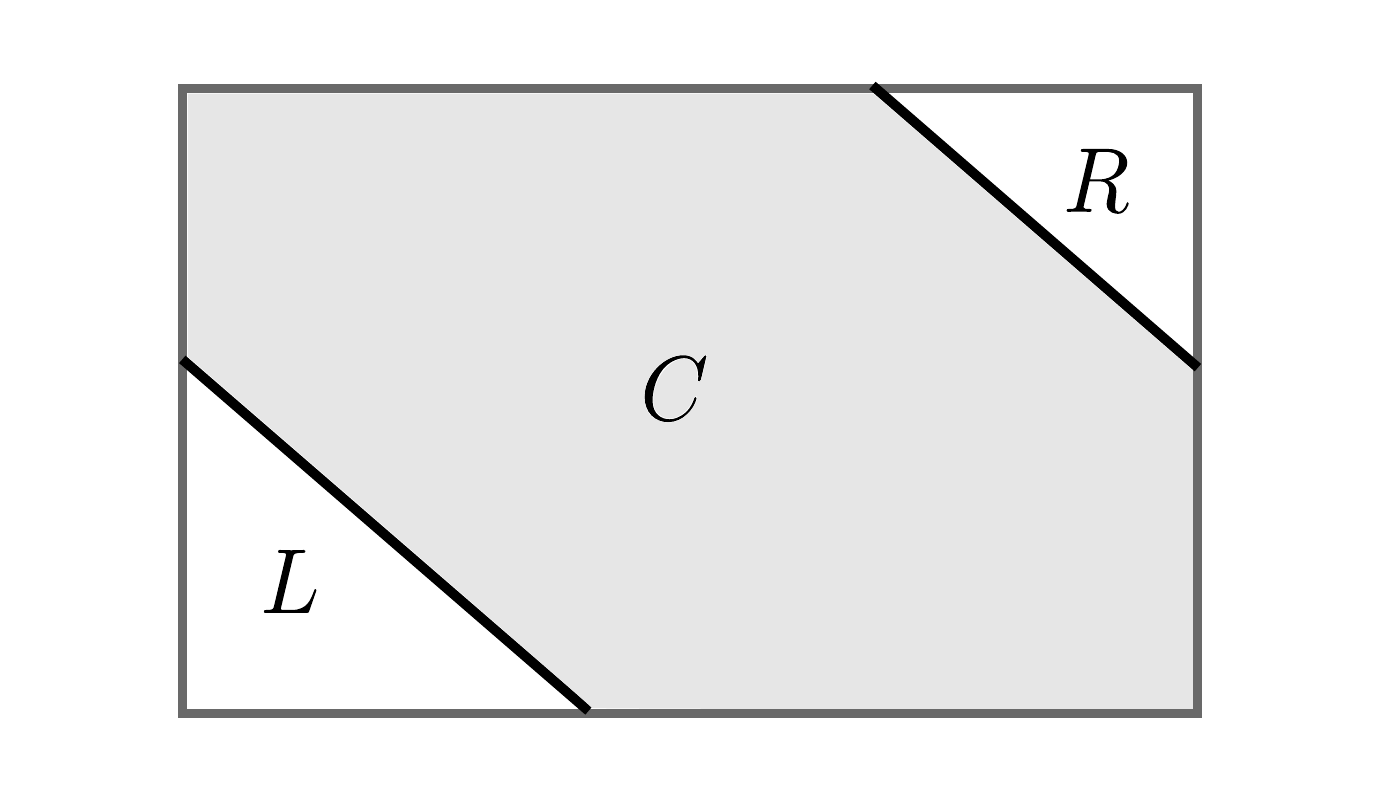}
	\end{subfigure}	
	\caption{Integration over a polygonal region $C$ whch is a subset of the cell  $L\cup C \cup R$. For the integration over $C$, we subtract the integrals over $L$ and $R$ from that over $L\cup C \cup R$.}\label{parallelogram intergration}
\end{figure}

\section*{Acknowledgement}

S. Y. Cho, S. Boscarino and G. Russo would like to thank the Italian Ministry of Instruction, University and Research (MIUR) to support this research with funds coming from PRIN Project 2017 (No. 2017KKJP4X entitled “Innovative numerical methods for evolutionary partial differential equations and applications”). S. Boscarino and G. Russo are members of the INdAM Research group GNCS.

\bibliographystyle{amsplain}
	\bibliography{references}
	
\end{document}